\ifx\shlhetal\undefinedcontrolsequence\let\shlhetal\relax\fi


\def\today{\ifcase\month\or January\or February\or March\or
April\or
    May\or June\or July\or August\or September\or October\or
November
    \or December\fi\space\number\day, \number\year}
\magnification=1200
\font\bigbf=cmb10 scaled \magstep2
\def\nz{\hfil\break\noindent}
\def\nl{\nz}
\def\lz{\vskip6mm \noindent}
\def\hz{\vskip7.2pt \noindent}
\def\USP{{\rm USP}}

\def\tcf{{\rm tcf\,}}
\def\l{{\lambda}}
\def\ord{{\rm \ord\,}}
\def\lub{{\rm lub}}
\def\OP{{\rm OP\,}}

\font\msxm=msxm10 scaled 1200
\font\msym=msym10 scaled 1200

\def\pp{{\rm pp}}

\def\rk{\rangle}
\def\lk{\langle}
\def\lg{{\rm \ell g}}
\def\ap{{\rm ap}}

\def\qed{\hfill\hbox{\msxm\char'003}}
\def\Vdash{\kern2pt{\vrule depth 1pt height .8 em width .4pt}\kern
.1pt\vdash}
\def\notVdash{\;\>\rlap/\kern-5pt\Vdash}
\def\notvdash{\;\>\rlap/\kern-6pt\vdash}
\def\uhr{\hbox{\msxm\char'026}}
\def\name#1{\lower 1em\rlap{\char'176}#1}
\def\lq{{\rm``}}

\def\and{{\, \&\,}}
\font\msym=msym10
\textfont9=\msym
\font\msyms=msym7
\scriptfont9=\msyms
\def\bet{{\fam9 i}}

\def\bbq{{\fam9 Q}}

\def\Min{{\rm Min\,}}
\font\fatone=cmssbx10

\def\bK{\hbox{\fatone K}}

\def\eqdf{\buildrel\rm def\over =}
\def\d{\delta}
\def\beth{\bet}
\def\k{\kappa}
\def\l{\lambda}
\def\a{\alpha}
\def\w{\omega}
\def\({\big(}
\def\){\big)}
\def\SP{{\cal P}}
\def\SS{{\cal S}}

\def\CT{{\cal T}}
\def\Dom{{\rm Dom\,}}
\def\nacc{{\rm nacc\,}}
\def\acc{{\rm acc\,}}
\def\md{{\rm md}}
\def\vare{\varepsilon}
\def\id{{\rm id\,}}
\def\bd{{\rm bd}}

\def\cf{{\rm cf\,}}
\def\tlim{{\rm tlim\,}}
\def\Rang{{\rm Rang\,}}

\def\cov{{\rm cov\,}}

\def\otp{{\rm otp\,}}
\def\mod{{\rm mod\,}}
\def\ord{{\rm ord\,}}
\def\stat{{\rm stat\,}}

\def\univ{{\rm univ\,}}

\def\epsilon{\varepsilon}

\font\eufb=eufb10
\def\GB{\hbox{\eufb \char "42}}
\hfuzz3pt
\def\shelahtitle{\centerline{\bigbf The Universality Spectrum
: Consistency for more Classes}
\bigskip

\baselineskip24pt

\centerline{no. 457}
\bigskip
\centerline{Saharon
Shelah}
\medskip
\centerline{\sl  Institute of Mathematics, The Hebrew University,
Jerusalem, Israel}
\centerline{\sl Department of Mathematics, Rutgers University,
New Brunswick, N.J., U.S.A.}}

\def\BSF{{\rm Partially\; supported\; by\; the \;United\; States\;
Israel\;
binational\; science\; foundation}}
\def\feq{{\rm feq}}

\nl
 \shelahtitle
\footnote{}{$\BSF$, Publication No.~457,
\S\S2, 3 are some years old, 5.2,5.3 done in summer 92,
\S4 +5.1,  first version  done with \S\S2, 3 but written in sp.~92, 11.92
resp.. Revised 12/93.}
\def\univ{{\rm univ}}
\def\Univ{{\rm Univ}}
\nl
\centerline{Last revised, 1 August, 1993}
\bigskip
\noindent{\bf Abstract.} We deal with consistency results for
the existence of universal models in natural classes of models
(more exactly--a somewhat weaker version).
We apply a result on quite general family
to $T_{\feq}$ and to 
 the class of triangle-free graphs 
\bigskip
\noindent
{\bf \S0 Introduction:}
\medskip
The existence of universal structures, for a class of structures
 in a given cardinality is quite natural as witnessed by  having
arisen in many contexts. We had wanted here to peruse 
it in the 
general
context of model theory but almost all will interest a combinatorialist
who is just interested in the existence of universal linear order
or a triangle free graph.
For a first order theory (complete for simplicity) we look at the
universality spectrum
USP$_T=\{\lambda:\;T$ has a universal model in cardinal $\lambda\}$
(and variants). Classically we know that under GCH,
every $\lambda>\vert T\vert$ is in USP$_T$,
moreover
$2^{<\lambda}=\lambda>\vert T\vert\Rightarrow
\lambda\in \USP_T$
(i.e.--the existence of a saturated or special model, see e.g. [CK]).
Otherwise in general it is ``hard'' for a theory $T$ to have a
universal model
(at least when $T$ is unstable).
For consistency see [Sh100],
[Sh175],
[Sh 175a], Mekler [M] and parallel to this work Kojman-Shelah
[KjSh 456]  ;
on ZFC nonexistence results see
Kojman-Shelah [KjSh409],
[KjSh447],
[KjSh455].
We get ZFC non existence result
 (for $T^*_{{\rm feq}}$
under more restriction , essentially cases of failure of SCH
) in \S2,
more on linear orders (in \S3),
consistency of (somewhat weaker versions  of)
existence results abstractly
 (in \S4) derived consistency results
and apply them to  the class of models of
$T_\feq$  (an indexed family of independent equivalence relations)  
and  to the class of  triangle free graphs  (in \S5 ).
The general theorem in \S4  was  intended for treating all simple 
theories (in the sense of [Sh 93] , but this is  not included as it is 
probably too much model theory for the expected reader here
 (and for technical reasons).
\medskip
\noindent
{\bf \S1}
\nl
{\bf 1.1 Definition:}
 For a class $\bK=(\bK, \le_{{\bK}})$ of models
\nl
1) $\bK_\lambda=\{M\in \bK:\Vert M\Vert=\lambda\}$
\nl
2) $\univ (\lambda, \bK)=\Min \{ \vert {\cal P}\vert: {\cal P}$
a set of models from $\bK_\lambda$ such that for every
$N\in \bK_\lambda$ for some $N\in {\cal P}$,
$M$ can be $\le_{{\bK}} $-embedded into $N \}$.
\nl
3) $\Univ (\lambda, \bK)=\Min \{\Vert N\Vert:
N\in \bK$ , and every $M\in \bK_\lambda$ can be
$\le_{\bK}$-embedded into $N\}$.
\nl
4) If $\bK$ is the class of models of $T$, $T$ a complete theory, we
write $T$ instead $(\mod T, \prec )$
(i.e. the class of model of $T$ with elementary embeddings). If $\bK$ is the class of models
of $T$,
$T$ a universal theory, we write $T$ instead $(\mod(T), \subseteq)$.
\lz
{\bf 1.2  Claim:} 1)
$\univ(\lambda, \bK)=1$ iff
$\bK$ has a universal member 
of
cardinality $\lambda $.
\nl
2) Let $T$ be  first order complete,
$\vert T\vert\le \lambda$.
{\it Then} we have
$\univ(\lambda, T)\le\lambda$  implies  
$\univ(\lambda, \bK) =1$ and 
$\Univ(\lambda, T )  \le \univ(\lambda,  T)\le
\cf \( \SS_{\le  \lambda } ( \Univ (\lambda , T ) , \subseteq) \) =
\cov(\Univ(\lambda,  T), \lambda^+, \lambda^+, 2)$
(see [Sh-g] ; we can replace $T$ with  $\bK$
with
suitable
properties).
\lz
{\bf \S2 The universality Spectrum of $T_{{\rm feq}}$}

For $T_{{\rm feq}}$, a prime example for a  theory with the tree order 
property
(but not the strict order property), we prove there are limitations on
the universality spectrum;
it is meaningful when SCH fails.
\hz
{\bf 2.1 Definition:}
$T^*_{\feq}$ is the model completion of the following
theory,
$T_{\feq}$.   $T_{\feq}$ is defined
as follows:
\item{(a)} it has  predicates $P, Q$ (unary)
$E$ (three place,
written as $yE_xz\}$
\item{(b)}  the universe (of any model of $T$)
is the disjoint union of $P$ and $Q$
, each infinite
\item{(c)} $yE_xz\to P(x)\,\&\, Q(y)\,\&\, Q(z)$
\item{(d)} for any fixed $x\in P$,
$E_x$ is an equivalence relation on $Q$
with infinitely many equivalence classes
\item{(e)} if $n<\omega$,
$x_1, \dots, x_n\in P$ with no repetition and
$y_1, \dots, y_n\in Q$ {\it then} for some $y\in Q$,
$\bigwedge_{\ell=1}^n yE_{x_\ell} y_\ell$.
\nl
(Note: $T_{\feq}$ has elimination of quantifiers).
\lz
{\bf 2.2 Claim:}
Assume:
\item{(a)} $\theta<\mu<\lambda$
\item{(b)} $\cf\lambda=\lambda, \theta=\cf\theta=\cf\mu$,
 $\mu^+<\lambda$
\item{(c)} $\chi=:pp_{\Gamma(\theta)} (\mu)>\lambda
+ \vert i^* \vert$
\item{(d)} there is $\{(a_i, b_i):i<i^*\}$,
 $a_i\in [\lambda]^{<\mu}$,
$b_i\in [\lambda]^\theta$ 
and  $ \vert \{ b_i : i< i^* \}  \vert \le \l $
such that:
for every $f:\lambda\to \lambda$ for some $i$,
$f(b_i)\subseteq a_i$

\noindent
{\it then}
\nl
(1) $T_{\feq}$ has no universal model in $\lambda$.
\nl
(2) Moreover,
$\univ(\lambda, T_{\rm feq})\ge \chi=pp_{\Gamma(\theta)}(\mu)$ .
\hz
{\it Proof:}
Let $D$ be a $\theta$-complete filter on $\theta$,
$\lambda_i=\cf\lambda_i<\mu=\sum_{i<\k} \lambda_i$,
$\tlim_{D}\lambda_i=\mu,$
$\chi=:$  tcf($\prod_{i<\theta}\lambda_i/D)>i^*$
(and for (2),
$\tcf \( \prod_{i<\theta}\lambda_i/D  \)>
\univ(\lambda, T_{{\rm feq}}))$.
Also
let $\lk f_\alpha:\alpha<\chi\rk $
be $ <_D$-increasing cofinal in $\prod_{i<\theta}\lambda_i/D$.
Let $S=\{\delta<\lambda:\cf\delta=\theta,
\delta$ divisible by $\mu^{\omega+1}\}$.
Let $\bar C=\lk C_\delta:\delta\in S\rk$
be such that:
$C_\delta$ a club of $\delta$,
otp$(C_\delta)=\mu$ and
$[\alpha\in C_\delta\Rightarrow \alpha>0$ divisible by
$\mu^\omega$]
and $\emptyset \notin \id^a(\bar C)$
(i.e. for every club $E$ of $\lambda$ for
stationary many $\delta\in S\cap E$,
$C_\delta\subseteq E$) (exists-see [Sh 365,\S2]).

For (1), let $M^*$ be a candidate for being a universal model of $
T_{\feq}$
of cardinality $\lambda$ , for (2) let
$\langle M^*_\zeta:\zeta<\kappa\rangle$
exemplify $\kappa=:\univ (\lambda, T_{\rm feq});$
for (1) let
$\kappa=1$,
$M^*_0=M_0$.
Without loss of generality $\vert P^{M^*_\zeta}\vert=\vert
Q^{M^*_\zeta}\vert=\lambda$,
 $P^{M^*_\zeta}$ is the set of even ordinals $<\lambda$,
$Q^{M^*_\zeta}$ is the set of odd ordinals $<\lambda$.

For each $i<i^*$
and
$\delta\in S$ and
$z\in Q^{M^*_\zeta}$ let
 $a'_i=\{2\alpha:\alpha\in a_i\}$ and
$d[z, \delta, i, \zeta]=\{\alpha:\alpha\in \nacc C_\delta$
and for some $x\in a'_i$ there is $y<\alpha$,
such that $M^*_\zeta\models y E_x z$ but there is no
$y<\sup(\alpha\cap C_\delta)$ such that
$M^*_\zeta\models yE_xz\}$.
Clearly $d[z, \delta, i, \zeta]$ is a subset of $C_\delta$ of
cardinality $\le\vert a_i\vert<\mu$.

Define $g_{z, \delta, i, \zeta}\in \prod_{j<\theta}\lambda_j$ by:
if $\vert a_i\vert <\lambda_j,
\beta\in C_\delta, \;\otp (\beta\cap C_\delta)=\lambda_j$
then $g_{z, \delta, i, \zeta}(j)=\otp
(\vare\cap C_\delta)$
where $\vare\in C_\delta\cap \beta$ is
$\Min\{\vare:\vare\in C_\delta\cap \beta$,
 $\vare>\sup(d[z, \delta, i, \zeta]\cap \beta)\}$ and let
$g_{z, \delta, i,\zeta }(j)=0$
if $\vert a_i\vert\ge \lambda_j$.
 By the choice of $\langle f_\alpha:\alpha<\chi\rangle$ for some
$\gamma$ we have
$g_{z, \delta, i, \zeta}<_Df_{\gamma}$, let
$\gamma^*=\gamma^*[z, \delta, i, \zeta]$ be the first such $\gamma$.
As $\mu=\tlim_D \lambda_i$
clearly  $\gamma^*[z, \delta, i, \zeta]$ is the first
$\gamma<\chi$ such that for the
$D$-majority of $i<\theta$,
$\bigwedge_{\alpha\in d[z, \delta, i, \zeta]}\otp
(\alpha\cap C_\delta)\notin[f_\gamma(i), \lambda_i)$;
clearly it is well defined.
Wlog 
$\{ b_i : i < i^* \}  = \{ b_i : i < i^* \}\cap \lambda $

As $\chi>\lambda   
+\k + \vert i^* \vert$,
 there is $\gamma (*) <\chi$ such that:
$z\in Q^{M^*_\zeta}$,
$\delta\in S$,
 $i<i^*$, $\zeta<\kappa \Rightarrow \gamma^*[z, \delta, i,
\zeta]<\gamma (*)$.
Now we can define by induction on $\alpha<\lambda$,
$N_\alpha, \gamma_\alpha$ such that:

\item{(i)} $N_\alpha $ is a model of $T^*_{\feq}$ with universe
$\gamma_\alpha= \mu(1+\alpha)$,

\item{(ii)}
all $x\in P^{N_\alpha}$
are even,
all $y\in Q^{N_\alpha}$
are odd

\item{(iii)} $N_\alpha$ increasing continuous,
$P^{N_\alpha}\not= P^{N_{\alpha+1}}$

\item{(iv)} for any $x\in P^{N_\alpha}$
there is a $y=y_{x, \alpha}\in Q^{N_{\alpha+1}}\setminus 
Q^{N_\alpha}$ 
such that
$\neg(\exists z\in Q^{N_\alpha})[z E_xy]$,
\item{(v)} if
$\alpha\in S, i<\alpha\cap
 i^* \cap \l $
and  $b'_i  \subseteq \Min ( C_\alpha  )$
then there is a $z^i_\alpha\in
Q^{N_{\alpha+1}}\setminus
Q^{N_\alpha}$ such that
$\Rang f_{ \gamma (*) }=\{\otp(y\cap C_\alpha ):$ for some $x\in 
b'_i$
, $y$ is
minimal such that $yE_xz^i_\alpha\}$ where
$b'_i\eqdf\{2\alpha:\alpha\in b_i\}$.

\noindent  [For  carrying out this let
$ d_{\alpha , i } =^{df} \{  \beta \in C_\alpha  : \otp (
C_\alpha  \cap \beta ) = 
 ( f_{\gamma (*) } (j) +1 ) $ for some $ j<\theta \} $,
so $d_{\alpha, i}\subseteq \nacc(C_\alpha)$, now
 choose distinct
$x_{\alpha  , i, \beta } \in b'_i$  for $\beta \in d_{\alpha , i } $ .
Next choose
$y_{ \alpha , i , \beta } \in  \beta \setminus \sup ( C_\alpha \cap 
\beta )$
such that it is as in clause (iv)  for
$x_{\alpha  , i , \beta} $   and  $z^i_\alpha E_{x_{\alpha ,
i,\beta  }} y_{\alpha , i , \beta }  $ .]

If $\zeta<\kappa$ and $f$ is
 an embedding of $N=\bigcup_{\alpha<\lambda}N_\alpha$ into
$M^*_\zeta$, for some $i$ we have $f(b'_i)\subseteq a'_i$
as we can define $f':\lambda\to \lambda$ by
$f(2\alpha)=2 f'(\alpha)$, well defined as $f$ maps
$P^N$ into $P^{M^*_\zeta}$.
Let $i_1< \lambda$ be such that $b_{i_1}=b_i$.
Let
$E=\{\delta<\lambda:(M^*_\zeta\uhr\delta, N\uhr\delta, f)\prec
(M^*_\zeta, N, f)$ and
$\delta>i_1\}$,
clearly it is a club of $\lambda$ hence,
by
the choice of $\bar C$,  for some $\delta\in S$ we have $C_\delta\subseteq
E$.
Let $z\eqdf f(z^{i_1}_\alpha)$,
so
$d[z, \delta, i_1, \zeta]$ is well defined.
For each $j<\theta$ there are $\beta_0<\beta_1$
in $C_\delta$ such that
 $\otp(C_\delta\cap \beta_0)=f_{\gamma(*)}(j)$,
$\otp(C_\delta\cap \beta_1)=f_{\gamma(*)}(j)+1$ and there is $y$
in $(\beta_1\setminus \beta_0)\cap Q^N$
and $x\in b'_{i_1}(=b'_i)$ such that
$yE_x z^{i_1}_\alpha$,
$y$ minimal for those $z^\alpha_i,x$.
So
$x^*=f(x)\in a'_i\subseteq M^*_\zeta$,
$f(z^\alpha_i)\in M^*_\zeta$,
and
letting $y^*=f(y)$ we have
$y^*<\beta_1$,
and
$y^*E_{x^*}f(z^i_\alpha)$.
Is there $y^*_1<\beta_0$ with those properties?
if so
$f(y)E_{f(x)}y^*_1$,
$(M^*_\zeta, N, f)\models (\exists
t)[t\in Q^{M^*_\zeta}\,\&\,f(t)E^{M^*_\zeta}_{x^*}y^*_1]$ so as
$x<\beta_0
, y^*_1<\beta_0\in E_1$ there is such $t<\beta_0$,
as $E$ is an
equivalence relation
$f(t)E_{x^*}f(z^i_\alpha)$. Now as $f$ is an embedding
$tE^N_x z^i_\alpha$,
contradicting the choice of $y$.
So $ y^*=f(y)$ witness
$\beta_1\in d[z, \delta, i,\zeta)$ hence $\otp(\beta_1\cap
C_\delta)\le g_{f(z^{i_1}_\alpha), \delta, i_1, \zeta}(j)$

We easily
get a
contradiction.
\nl
\null\hfill \qed$_{2.2}$
\lz
{\bf 2.3 Claim:} 1) In 2.2 we can replace clauses (c), (d) by
(c)$^+$, (d)$^-$ below and the conclusions still hold.

\item{(c)$^+$} $\chi=pp_D(\mu)>\vert i^*\vert +\lambda$,
  $D$ a filter on $\theta$, or at least for some $\langle
\lambda_i:i<\theta\rangle$,
$\lambda_i=\cf\lambda_i<\mu={\rm tlim}_J\langle \mu_i:i<\theta\rangle$ and
$\prod_{i<\theta}\lambda_i / D$ is $\chi$-directed.

\item{(d)$^-$}
$\{(a_i, b_i):i<
 i^*  \} $ , 
$a_i\in [\lambda]^{<\mu}$, 
$i^* \le \l $ or at least $\{ b_i : i< i^* \} $ has cardinality $\le \l
$,
$b_i=\{\alpha_{i,\zeta}:\zeta<\theta\}$ and
for every
$f:\lambda\to \lambda$ for some $i$ we have
$\{\zeta<\theta: f(\alpha_{i, \zeta})\in a_i\}\not=\emptyset\ \mod D$.

\noindent
2) Above we can weaken in (c)$^+$ the demand
$\lq \prod_{i<\theta}\lambda_i/D$ is $\chi$-directed'' by
\nl
$\lq \cf(\prod_{i<\theta}\lambda_i/D)\ge \chi$''
if in clause (d)$^-$ we strengthen $\lq \models \emptyset \mod D$''
to ``$\in D$''.
\hz
Also similarly we can prove
\nz
{\bf 2.4 Claim:} Assume
\nl
(a) $\theta<\mu<\lambda\le\lambda^*$
\nl
(b) $\theta=\cf(\theta)=\cf(\mu)=\cf\lambda$,
$\mu^+<\lambda$
\nl
(c) $\pp_D(\mu)>\vert i^*\vert+\cov(\lambda, \mu, \theta^+,
\sigma)$, $D$ is $\sigma$-complete.
\nl
(d)$^-$ $\{(a_i, b_i):i<
 i^*  \} $ , 
$a_i\in [\lambda]^{<\mu}$, 
$i^* \le \l $ or at least $\{ b_i : i< i^* \} $ has cardinality $\le \l
$,
\nl\phantom{(d)$^-$}
$b_i=\{\alpha_{i,\zeta}:\zeta<\theta\}$ and
for every
$f:\lambda\to \lambda$ for some $i$ we have
$\{\zeta<\theta: f(\alpha_{i, \zeta})\in a_i\}\not=$
\nl\phantom{(d)$^-$}
$\emptyset\ \mod D$.
\nl
{\it Then}
\nl
(1) $T_{{\rm feq}}$ has no universal model in $\lambda$
\nl
(2) moreover $\univ(\lambda, T_{{\rm feq}})\ge\pp_D(\mu)$
\lz
{\bf 2.5 Remark:}
1) When  does (d) of 2.2 hold?;
it is a condition on $\l > \mu > \theta$,
assuming
 for simplicity $ \theta > \aleph_0$, $i^*=\lambda$)
 e.g. it holds (even with $\bigwedge b_i=b_0$) if:
\nl
\itemitem{$(*)_1$}
for some cardinal $\kappa$ we have
$\kappa^\theta\le \lambda$,
$\kappa=\cf\kappa$,
 $\cov (\lambda, \kappa^+, \kappa^+, \kappa)\le \lambda$.
\nl

\item{2)} As for condition $(d)^-$ from claim 2.3, if $D$ is the
 filter of co-bounded subsets of $\theta$ , it suffices to have
\itemitem{$(*)_2$} for some cardinal $\kappa$ 
we have
$\cov (\lambda , \mu , \kappa^+ , \kappa ) \le \lambda $ ,
or equivalently,
 $\sigma \in [\mu , \lambda )$ and $\cf (\sigma ) = \kappa$
imply 
$\pp_{\Gamma(\kappa )} ( \sigma ) \le \lambda $.
\item{3)} So if  $ \theta = \cf ( \mu ) < \beth_\omega ( \theta )
\le \mu < \mu^+ < \lambda = \cf (\lambda ) <
\pp_{\Gamma ( \theta )} ( \mu ) $
then by [Sh  460]  condition $(*)_1$ holds 
for some  $ \k < \beth_\omega (\theta ) $
\item{4)} Why have we require $\theta > \aleph_0$ ? as then 
by [ Sh-g, Ch. II , 5.4]  we can describe the instances of $\cov$ by 
instances of $\pp$ ; now even without this restriction this usually 
holds (see there) and possibly it always hold ; alternatively , we can 
repeat the proof of 2.2 using $\cov$
\item{5)} The parallel of 2.3(2) for 2.4 can be easily stated.
\lz
{\bf 2.6 Conclusion:}
If $\theta=\cf\mu$;
$\bet_\omega(\theta)\le\mu$,
$\mu^+<\lambda=\cf\lambda<\pp^+_{\Gamma(\theta)}(\mu)$ then $\univ
(\lambda, T_{\feq})\ge \pp_{\Gamma(\theta)}(\mu)$.
\hz
{\it Proof:}
The next step is:
\lz
{\bf 2.7 Question:}
Let $T$ be f.o. with the tree property without the strict order
property;
\nl
(see [Sh-c])
does 2.2 hold for it?
\lz
{\bf \S 3 A consequence of the existence of a universal linear order.}

This section continues, most directly, [KjSh 409].
\lz
{\bf 3.1 Claim:}
Assume

(a)$_\lambda$ $\kappa<\lambda\le 2^\kappa$ and
$2^{<\lambda}   \le  \lambda^+ < 2^\lambda$,
 $\lambda$ is regular.

(b)$_\lambda$ in $\mu=\lambda^+$ there is a universal linear order
\nz
{\it then}
\item{$\otimes_{\lambda , \mu }$}
there are $f_\alpha:\lambda\to \lambda$
(for $\alpha<\mu $) such that:
\nz
$(*)_{\lambda , \mu }$ for no $f:\lambda\to \lambda$ do we have
$\bigwedge_{\alpha<\mu } f_\alpha \not=_{J^{\bd}_\lambda}f$.
\hz
{\it Proof:}
Assume $\otimes_{\lambda , \mu }$ fails.
We use $\kappa$-tuples of elements to compute invariants.
Note that
$2^\kappa \le 2^{< \lambda} \le \lambda^+$
hence
$2^\kappa \in \{ \lambda , \lambda^+ \}$
hence $ (\lambda^+ )^\kappa = \lambda^+ $ .
Let $\lk\bar x^\vare:\vare<\lambda^+\rk$
list ${}^\kappa(\lambda^+)$.
Let $\lk \eta_\alpha:\alpha<\lambda\rk$
list
$\lambda$ distinct members of ${}^\kappa 2$
(not necessarily all of them ).
 Note that as $2^{<\lambda}\le \lambda^+$
there is
a stationary 
$S\in I[\lambda]$, $S\subseteq \{\delta<\lambda^+:
\cf(\delta)=\lambda\}$ (see [Sh 365 , \S2] for the definition of $I[S]$).

As $S\in I[\lambda]$ by [Sh365, \S2] there is
$\bar C=\lk C_\delta:\delta\in S\rk$ an $S$-club system such that
$\emptyset \notin\id_p(\bar C)$,
 $\otp C_\delta=\lambda$ and 
\nl
$\oplus$ for each $\alpha<\lambda$ we have
$\vert \{C_\delta\cap \alpha:\alpha \in\nacc
C_\delta\}\vert\le\lambda
$.

Let $M^*$ be a candidate for being a universal model of $T_{\ord}$
of cardinality  $\lambda^+$, wlog with universe  $\lambda^+$ .

For every linear order $M$ with universe $\lambda^+$, for every
$\bar
x \in
{}^\kappa M$
(a $\kappa$-tuple of members of $M$)
and $\delta\in S$, we define
a (possibly partial) function $g=g^{\bar x}_{M,\delta}:
\nacc C_\delta\to \lambda$ as follows:
\item{$(*)_0$}
for $\alpha\in\nacc C_\delta$,
$g(\alpha)=\beta$  {\it iff}  for every 
$\zeta<\kappa$ we have :
\nl
$\eta_\beta(\zeta)=1\iff
(\forall \gamma<\alpha)(\exists \gamma'<\sup(\alpha\cap C_\delta)$)
[$\gamma<_M x_\zeta\Rightarrow \gamma<_M\gamma'<_M x_\zeta]$.

Clearly $g^{\bar x}_{M, \d } ( \a ) $ can have at most one value . 
We call $(\delta, \bar x)$ {\sl good} in $M$ if for every $\alpha\in 
\nacc
C_\delta$
there is $\vare<\delta$ such that :
$ \bar x^\vare$,
$\bar x$ realize the same $<_M$-Dedekind cut 
over
$\{i:i<\sup(\alpha\cap C_\delta) \} $
(necessary if $2^{<\lambda}=\lambda^+$).
(The 
 meaning is that for every $\zeta<\kappa$,
$x^\epsilon_\zeta$, $x_\zeta$ realize the same $<_M $-Dedekind cut  
over  $\{i:i<\sup(\alpha\cap C_\delta) \} $  ).

Let $h_\delta:\lambda\to \nacc C_\delta$ be:
$h(i)$ is the $(i+1)$-th member of $C_\delta$.
We are assuming
``$\otimes_{\lambda , \mu }$ fails'',
so $\{g^{\bar x}_{M^*,\delta}\circ h_\delta:\bar x\in {}^\kappa 2,$
$\delta\in S\}$
cannot exemplify it.
So 
we can find $h^*_{M^*} :\lambda\to \lambda$ such that:
\item {$\otimes$} if $\bar x\in {}^{\kappa}(M^*) $,
$\delta\in S$ is
$(\delta, \bar x)$ good in $M^*$
then $(g^{\bar x}_{M^* , \delta}\circ h_\delta)\in {}^\lambda\lambda$
satisfies
$h^*\not=_{J_\lambda^\bd}(g^{\bar x}_{M^* , \delta}\circ h_\delta)$.

Let  $h^* = h^*_{M^*} $ ;
let $g_\delta:\nacc C_\delta\to \lambda$ be $h^*\circ (h^{-1}_\delta):
\nacc C_\delta\to \lambda$.
We now as in
 [KjSh 409, x.x?????] (using $S\in I[\lambda]$  i.e.  $\oplus$ )
construct a
linear order
$N=M^{h^*}$ with universe  $\lambda^+ ,$
$N=\bigcup_{\alpha<\lambda}N_\alpha$,
$N_\alpha$ increasing continuous in $\alpha$ with universe
 an ordinal $<\lambda^+$ and for each $\delta\in
S$,
there is a sequence $\bar y^\delta=\lk
y^\delta_\zeta:\zeta<\kappa\rk$
of members of $N_{\delta+1}$ such that
\item{$(*)_1$}
if $\alpha\in \nacc C_\delta$,
$g_\delta(\alpha)=\beta$,
$\zeta<\kappa$ {\it then}
\nl
$\eta_\beta(\zeta)=1\Leftrightarrow
(\forall\gamma\in N_\alpha )(\exists \gamma'\in N_{\sup(\alpha\cap
C_\delta)})
[\gamma<_N  y^\delta_\zeta\Rightarrow
\gamma<_N\gamma'<_N y_\zeta^\delta]$.

Suppose $f:\lambda^+ \to \lambda^+ $ is an embedding of $N$
into $M^*$,
let $E=\{\delta < \lambda^+ :
N_\delta$ universe is $\delta$ and $\delta$ is closed under $f$,
$f^{-1}\}$.
Clearly $E$ is a club of $\lambda^+$,
hence for some $\delta\in S$ the set
$A=(\acc E)\cap (\nacc C_\delta)$ is unbounded in $\delta$
(so $\delta\in \acc\acc E)$.
Let
$\bar x =\lk x_\zeta:\zeta<\kappa\rk =:\lk
f(y^\delta_\zeta):\zeta<\kappa\rk$,
so we know
(similarly to [KjSh 409 \S3]????)
that for $\alpha\in A$ and
$\zeta<\kappa$ we have
${g^{\bar x}_{M^* , \delta}(\alpha)}(\zeta)=1\Leftrightarrow
\eta_{g_\delta(\alpha)}(\zeta)=1$.
Hence $\alpha\in A\Rightarrow g^{\bar x}_{M^* ,
\delta}(\alpha)=g_\delta(\alpha)\Rightarrow
(g^{\bar x}_{M^*, \delta}\circ h_\delta)(\otp(\alpha\cap C_\delta) - 1)=
h^*(\otp(\alpha\cap C_\delta) - 1 )$  contradicting the choice of $h^*$.
\hfill\qed$_{3.1}$
\lz
{\bf 3.1A Claim:}
1) In 3.1 if $\lambda$ is a successor cardinal then we can
get
\nl
 $\oplus^0_\lambda$ there are
$f_\alpha:\lambda\to \lambda$ for $
\alpha<\lambda^+$ such that
\item{$(*)^\lambda$} for every
$f\in {}^\lambda\lambda$ for some $\alpha<\lambda^+$ we have
$f_\alpha\not=_{D_\lambda}f$
(where $D_\lambda$ is the club filter on $\lambda$).
\hz
2) If we allow $\mu>\lambda^+$, clause (a) of 3.1 holds  and
$(b)^*$ below then $\otimes_{\lambda , \mu }$ of 3.1  holds  ;
similarly in 3.1A(1), where
\nl
$(b)^*$ 
$ \univ ( \lambda^+ , T_{\ord} )\le \mu $
\hz
{\it Proof:} 1) Use [Sh 413, 3.4].
\nl
2) The same proofs.

So from the existence of a universal linear order of cardinality
$\lambda^+$, where
$\lambda$  is as in 3.1+3.1A(1), we get $\oplus^\lambda$,
from this we get below a stronger guessing of clubs.
\hfill\qed$_{3.1A}$
\lz
{\bf 3.2 Claim:}
Assume $\lambda$ is regular uncountable
 , and
\item{$\otimes^1_\lambda$} there are $f_\zeta:\lambda\to \lambda$
for $\zeta<\lambda^+$ such that:
for every
$f:\lambda\to \lambda$
for some $\zeta$,
$\{\alpha < \lambda:
f_\zeta(\alpha)=f(\alpha)\}$ is stationary.
\par\noindent
1) Let $S_1\subseteq \{\delta<\lambda^+$:
$\cf(\delta)=\lambda\}$, $S_2\subseteq\lambda$ be
stationary,
 and $\delta\in S_1\Rightarrow \delta=\sup(\delta\cap S_2)$.
We can find $\bar C=\lk C^\zeta_\delta:\delta\in S_1,
\zeta<\lambda^+\rk$,
such that :

(a) $C^\zeta_\delta $ is a club of $\delta$ of order type $\lambda$.

(b) $\nacc C^\zeta_\delta\subseteq S_2$.

(c) for every club $E$ of $\lambda^+$,
for stationarily many
$\delta\in S_1$,
for some $\zeta<\lambda^+$,
\nl\phantom{mm (c)} $\delta=\sup\Big\{\alpha:
\alpha\in \nacc C^\zeta_\delta$
and $\sup(\alpha \cap C^\zeta_\delta)\in \nacc C^\zeta_\delta$,
$\otp (\alpha\cap C_\delta)$ is
\nl\phantom{mm (c)
$\delta=\sup\Big\{\alpha:$}
even and
$\{\alpha, \sup(\alpha\cap C^\zeta_\delta)\}\subseteq E\Big\}$.

\noindent
2) Let $\lambda=\lambda^{<\lambda} $ 
and 
$ S \subseteq \{\delta<\lambda^+:\cf \delta=\lambda\}$
 stationary.
We can find $\bar C=\lk C^\zeta_\delta:\delta\in S,
\zeta<\lambda^+\rk$
such that
\itemitem{(a)} $  C^\zeta_\delta $ is a club of $\delta$ of order type
$\lambda$.
\itemitem{(b)} for every club $E$ of $\lambda^+$
for stationary many $\delta\in S$,
for some $\zeta<\lambda^+$,
for every $ \xi < \l $ we have
$E$ contains arbitrarily large 
(below $\l$ )
intervals of    $  C_\delta  $ of length
$\xi$
\par\noindent
3) If $\lambda$ is a successor cardinal then
we can get
(2) even if we omit ``$\lambda=\lambda^{<\lambda}$'' and weaken in
$\otimes^1_\lambda$ ,
``$f_\zeta(\alpha)=f(\alpha)$'' to
``$f_\zeta(\alpha)\ge f(\alpha)$''.
\par\noindent
4) In part(2), if $S_2=\lambda^+$ we can omit
``$\lambda=\lambda^{<\lambda}$'' if we restrict ourselves in (b) to $\xi$
a regular cardinal.
\noindent
\hz
{\bf  3.2A  Remark }
1) We can in 3.2(3) get the conclusion of 3.2(2) too if we fix $\xi$
\nl
2) We can replace in the assumptions and conclusions   ,
$\lambda^+$  by  $\mu$ is in 3.1A(2).
\hz
{\it Proof:}
1) Let $\lk C_\delta:\delta\in S_1\rk$
be such that:
$C_\delta$ a club of $\delta$,
$\otp C_\delta=\lambda$
and  $\nacc ( C_\delta ) \subseteq S_2$.
 If $\alpha<\beta<\lambda^+$,
$S_2\cap (\alpha, \beta)$ has at least two elements then
let $\lk(\beta^\vare_{\alpha, \beta}, \gamma^\vare_{\alpha,
\beta}):\vare<\lambda\rk$
 list all increasing pairs from $(S_2\cap \beta\setminus\alpha)$
(maybe with repetitions).
Let $\lk f_\zeta:\zeta<\lambda^+\rk$
exemplify $\otimes^1_\lambda$.
Let $C_\delta=\{\alpha_{\delta, \vare}:\vare<\lambda\}$
(increasing). Let $e=e^\zeta_\delta\subseteq \lambda$ be a club of
$\lambda$
such that:
if $i<j$ are from $e$ then
$\gamma^{f_\zeta(i)}_{\alpha_{\delta, i},  \delta } <
\alpha_{\delta, j}$.
Now for $\delta\in S_1$,
$\zeta<\lambda^+$,
we let:
$C^\zeta_\delta=
\{\alpha_{\delta, \vare} , \beta^{f_\zeta(\vare)}_{\alpha_{\delta, \vare},
 \delta},
\;\gamma^{f_{\zeta}(\vare)}_{\alpha_{\delta,\vare}, \delta
}:
\vare\in e^\zeta_\delta \}$.

Clearly $C^\zeta_\delta$ is a club of $\delta$ of order type
$\lambda$.
 Now if $E$ is a club
of $\lambda^+$, then $E\cap S_2$
 is a stationary subset of $\lambda^+$ so for some $\delta\in S_1$,
$\delta=\sup(E\cap S_2)$
 and define $g:\lambda\to \lambda$ by:
 $\beta^{g(\vare)}_{\alpha_{\delta,\vare} ,\delta }$,
 $\gamma^{g(\vare)}_{\alpha_{\delta, \vare}+1 ,\delta }$
are the
 first and second members of $(E\cap S_2)
\setminus  (\alpha_{\delta, \vare} ,
\delta ) \}$.
By the choice of $\lk f_\zeta:\zeta<\lambda^+\rk$ for some
$\zeta<\lambda^+$,
$(\exists^\stat \vare)(g(\vare) = f_\zeta (\vare))$.
So $C^\zeta_\delta$ is as required .
\nz
2) Similar proof (and we shall not use it).
\nz
3) In the proof of (1) for $\alpha<\lambda$
let $h(\alpha , -):
\lambda^-{\buildrel\rm onto\over\to} \alpha$.
We do the construction for each $\tau<\lambda^-$.
The demand on $e = e^\zeta_\delta$ is changed to:
if for $i<j$ are from $e$,
then $\gamma^{h(f_\zeta(\alpha), \tau)}_{\alpha_{\delta, i},
\delta} < \alpha_{\delta, j}$,
and $C^\zeta_\delta$ is changed accordingly.
For some $\tau<\lambda$ we succeed
(really this  version  of $\otimes^1_\l$ implies the original version .)
\nz
4) By the proof above we can get $C^\zeta_\delta$ such that:
for every regular $\xi<\lambda$ and club $E$ of $\lambda^+$ for
stationarily many
$\delta\in S_1$,
for unboundedly many $\alpha\in \nacc C^\zeta_\delta$,
we have:
$\alpha\in E$,
$\cf(\alpha)=\xi$.
Then we ``correct'' as usual (see[Sh365 \S2]).
\hfill\qed$_{3.2}$.
\lz
{\bf 3.3 Claim:}
Assume:

\item{(a)} $\lambda$ regular, $S\subseteq \lambda$ stationary,
$\lambda^\kappa=\lambda$.

\item{(b)} $\bar C=\lk C_\delta:\delta\in S\rk$,
 $C_\delta$ a club of $\delta$.

\item{(c)} $\bar \SP=\lk \SP_\delta:\delta\in S\rk,
\SP_\delta\subseteq \SP(\nacc(C_\delta))$ is closed upward.

\item{(d)} for every club $E$ of $\lambda$ for some $\delta$,
$E\cap \nacc C_\delta\in \SP_\delta$

\item{(e)} $\kappa<\lambda$,
$T_\delta=\bigcup\{T_{\delta,\beta, \gamma}:\beta\le \gamma
, \{ \beta , \gamma \} \subseteq  \nacc
C_\delta\}$ ,
for $\beta<\gamma\in \nacc C_\delta$,
$T_{\delta, \beta, \gamma}\subseteq
{}^{\gamma\cap \nacc C_\delta}({}^\kappa \beta)$,
$\vert T_{\delta, \beta, \gamma}\vert \le\lambda$,
and even
for each $\gamma$ the set $\bigcup \{T_{\delta, \beta,
\gamma}:\gamma\in\nacc C_\delta,\;
\beta\in \gamma\cap\nacc C_\delta\}$
has cardinality $\le\lambda$.

\item{(f)} If $A\in\SP_\delta$,
for $\zeta<\lambda^+$ we have
$f_\zeta\in {}^{\nacc C_\delta}({}^\kappa\delta)$ and
$\Big[\beta<\gamma$ are from $A\Rightarrow f_\zeta\uhr \beta \in
T_{\delta,
\beta, \gamma}\Big]$.
{\it Then } for some
$f^*
\in {}^{\nacc C_\delta}({}^\kappa \delta)$ we have
[$\beta<\gamma$ from $A\Rightarrow f^*\uhr \beta \in T_{\delta,
\beta,
\gamma}]$
and for every $\zeta<\lambda^+$,
$\{\beta\in A:f_\zeta(\beta)=f^*(\beta)\}\notin \SP_\delta$.
\hz
{\it Then} there is no universal
linear order of cardinality $\lambda^+$.
\hz
{\it Proof:} Similar to the previous one.

\hz
{\bf 3.5 Conclusion:}
If $2^\lambda > \lambda^+,\lambda=\cf \lambda>\aleph_0$,
$\bar C=\lk C_\delta:\delta\in S\rk$,
$S\subseteq \{\delta<\lambda^+:\cf \delta=\lambda \}$ stationary,
$\lambda^+\notin \id^a(\bar C)$ and for each
$\alpha$ we have $\vert\{C_\delta\cap \alpha:\alpha\in \nacc
C_\delta\}\vert\le \lambda$
{\it then}
\item{(a)}  there is no universal linear order in $\lambda^+$
\item{(b)} moreover , $\univ( \lambda^+ , T_\ord ) \ge 2^\lambda )$.
\lz
{\bf 3.6 Discussion:}
(1) The condition $\otimes_\lambda$ from 3.1 holds in the models (of
ZFC) constructed in [Sh 100, \S4]
where
$\lambda=\aleph_0$,
$2^{\aleph_0}=\aleph_2$ and there is a non meager subset of
${}^\omega 2$
of cardinality $\aleph_1$.
\nz
(2) It is clear from
3.5
that the existence of a universal graph in
$\mu$  does
not imply the existence of a universal linear order in $\mu$
every for $\mu=\lambda^+$,
$\lambda=\lambda^{<\lambda}$:
as
by [Sh 175], [Sh 175a],
if $V\models $ GCH,
$\lambda=\lambda^{<\lambda}$,
$\bar C=\lk C_\delta:\delta<\lambda^+, \cf \delta=\lambda\rk$
guesses clubs,
 for some $\lambda^+$-c.c. forcing notion
$P$ we have
 $V^P$ $\models_P$ ``there is a
universal graph
in $\lambda^+$''. But in $V^P$
the property of $\bar C$, guessing clubs,
is preserved and it shows that there is no universal linear order.
\nz
(3) We can look at this
from another point of view:
\item{(a)} Considering the following three proofs of consistency 
results  on 
the
existence
of universal structures:
[Sh 100, \S4]
(universal linear order in $\aleph_1$),
[Sh 175, \S1]
(universal graphs in $\lambda^+$,
$\lambda=\lambda^{<\lambda}$ and [Sh 175a]
(universal graphs in other cardinals),
the first result cannot be gotten by the other two proofs.
\item{(b)} For theories with the strict order property it is ``harder''
to have universal models than for simple theories (see [Sh93])
as the results of [Sh500, \S1] on simple theories fail for the theory of linear
order  (by 3.5) and even all (f.o.) theories with the strict order
property (as in [KjSh 409, x.x] )

\noindent(4) Concerning 3.5(b) ,  note that (for any complete first 
order $T$ )
we have $\Univ (\mu , T) \le 2^{< \mu } $ hence
$\cf \( \SS_{\le \mu } ( 2^{<\mu} ) , \subseteq \) \ge  \univ (\mu , T )$
so under reasonable hypotheses we get in 3.5(b)  equality (i.e.,  $\mu
=
\lambda^+ $).
\lz
{\bf \S4 Toward the consistency for simple theories}

The aim of this proof was originally to deal with the universality
spectrum of simple countable theories and as a first
approximation to 
characterize
$\{\lambda:\univ(\lambda^+ , T)\le\lambda^{++}< 2^\lambda\}$,
but we shall do it more generally and have more consequences. On 
simple theories see [Sh 93].
The reader may well read the ``smooth'' version, i.e. add in Definition
4.1,
the $(<\lambda)$-smoothness from 4.2(4), (5), and so we can omit clauses
(e)$(\beta),(\gamma), (\delta)+(1)$ from Definition 4.1. He can also
assume in 4.1 that $\tau_i=\tau_0$.
\lz
{\bf 4.0 Notation:}
(1) For a set $u\in {\cal S}_{<\lambda}(\lambda^+)=:\{
u\subseteq \lambda^+: \vert u\vert<\lambda\}$
let $\sup_\lambda(u)=\{\alpha +\lambda
:  \alpha\in u\}$ 
also  let
$S^{\lambda^+}_\lambda =\{\delta<\lambda^+:\cf\delta=\lambda\}$
\nl
(2) If $u_1, u_2\in \SS_{<\lambda}(\lambda^+),\; h:u_1\to u_2$ is {\it
legal} if it is one to one,
onto,
and there is a unique $h^+$ such that:
 $h^V$ is one to one order preserving
from
$\sup_\lambda(u_1)$
onto $\sup_\lambda(u_2)$ and for $\alpha\in u_1$,
$h^+(\alpha +\lambda )=h(\alpha)+\lambda$.
\nl
(3) We say that $h$ is {\it lawful} if in addition $h^+$ is the identity
 . We sometimes use
``legal" and  ``lawful" for functions $h : u_1  \to 
u_2$     when  $ u_i \subseteq \l^+ , \vert u_i \vert \ge \l $ .
\nl
(4) Wide $\lambda^+$-trees $\CT = ( \CT , < )$ are here-just subsets of
${}^{\lambda^+ >} (\lambda^+ )$  of cardinality $\le \lambda^+$
closed under initial segments with the order being initial segment.
A {\it branch} is a maximal linearly ordered subset, a {\it 
$\lambda^+$-branch}
is one of order type $\lambda^+$.
(So the trees are automatically normal).
\lz
{\bf 4.1 Definition:} $K_{{\rm ap}}=(K_{{\rm ap}},
\le_{K_{{\rm ap}}} ) $
is a $\lambda$-approximation family,
{\it if} for some sequence 
$ \bar \tau $
  ($= \lk \tau_i  : i< \lambda^+ \rk $ 
of vocabularies ,
$\vert \tau_i \vert\le\lambda$ , 
$\tau_i$ increasing with $i$,
$M\uhr i$
means $(M\uhr \tau_i)\uhr i$;
$\tau_i$ can have relations and functions with infinite
arity but $ < \lambda$
 (you may concentrate on the
case  $ \tau_i = \tau $ 
 for all $i<\l$  ) the following hold :

\item{(a)} $K_{{\rm ap}}$ is a set of $\tau$-model with a partial order
 $\le=\le_{K_{{\rm ap}}}$ (or $\mu$ is a $\tau_{\sup(M)}$-model).

\item{(b)} if $M\in K_{{\rm ap}}
$ then $\vert M\vert $ is a subset of $\lambda^+$ of cardinality
$<\lambda$ and  $ M \le_{K_{ap}}N \Rightarrow  M  \subseteq N $ .

\item{(c)} if $M\in K_{\ap},
\delta\in S^{\lambda^+}_\lambda$ then
$M\uhr \delta\in K_\ap$ and $M\uhr \delta\le_{K_{\ap}} M$;
also $M\uhr 0\in K_{\ap}$ (this is just to say we have the joint
embedding property).

\item{(d)} any $\le_{K_{\ap}}$-increasing chain in $K_\ap$ of length
$<\lambda$
has an upper bound.

\item{(e)} ($\alpha$) if $\delta\in S^{\lambda^+}_\lambda$,
$M_0=M_2\uhr \delta$,
$M_0\le_{K_{\ap}} M_1$,
$\vert M_1\vert \subseteq \delta$ {\it then}
$M_1$,
$M_2$ has a common $\le_{K_\ap}$-upper bound
$M_3$ ,  such that
$M_3\uhr \delta=M_1$
\nl
($\beta$) {\it if} we have
$M_{1, i}(i< i^*<\lambda)$,
$M_{1, i}\in K_{{\rm ap}}$ increasing with $i$ , 
$\vert M_{1, i}\vert\subseteq \delta_i\in S^{\lambda^+}_\lambda$
and 
$M_2\uhr\delta_i \le_{K_{\ap}} M_{1, i}$,
{\it then} there is a common upper bound $M_3$ to $\{M_2\}\cup
\{M_{1, i}:i<i^*\}$
\nl
($\gamma$) {\it if} we have $M_1 \in K_{\ap}$ ,
$M_{2,i} \in K_{\ap}$
for $i<i^* < \lambda$ increasing with $i$ ,
$\delta \in S^{\lambda^+}_\lambda$ ,
$M_{2,i} \uhr \delta \le M_1$
{\it then} there is a common $\le_{K_{\ap}}$-upper bound to
$\{ M_1\} \cup \{ M_{2,i} : i< i^* \}$
such that
$M_3\uhr \delta=M_1$ .
\nl
($\d $) {\it if} (i)
$\lk \d_i : i \le i^* \rk$
is a strictly increasing sequence of members of 
$S^{\lambda^+}_\lambda$ ,
\nl
\phantom{($\delta$) if}
(ii) we have
$M_{1, i}(i< i^*<\lambda)$,
$M_{1, i}\in K_{{\ap}}$ increasing with $i$ , 
\nl
\phantom{($\delta$) if}
(iii) $[ i(1) < i(2) \Rightarrow  M_{i(1)} = M_{i(2)} \uhr \d_{i(1)}  ] $
\phantom{(e)($\delta$) if }
\nl
\phantom{($\delta$) if}
(iv) $\vert M_{1, i}\vert\subseteq \delta_i$ 
\nl
\phantom{($\delta$) if}
(v) $M_{2,j} \in K_{\ap}$ for $j<j^*$
 has universe $\subseteq \d_{i^*}$ , 
and is $<_{K_{\ap}}$-increasing 
 in $j$
\nl
\phantom{($\delta$) if}
(vi) 
$M_{2,j} \uhr\delta_i \le_{K_{\ap}} M_{1, i}$,
\nl
{\bf then } there is a common upper bound $M_3$ to $\{M_{2,j} : j<j^* 
\}\cup
\{M_{1, i}:i<i^*\}$  such that for every  $i<i^*$
we have  $M_3 \uhr \d_i = M_{1, i } $

\item{(f)} For $\alpha <\lambda ^+$, 
$\{M\in K_{\ap}:\vert M\vert\subseteq \alpha\}$ has cardinality
$\le
\lambda$.

\item{(g)} We call $h:M_1\to M_2$ a {\it lawful}  ({\it legal})
$K_\ap$-isomorphism if $h$ is an isomorphism from $M_1$ onto
$M_2$ and $h$ is lawful (legal).
We demand:
\nl
\phantom{(g)}($\alpha$)
if $M_1\in K_{\ap}$,
$u_1=\vert M_1\vert$,
$u_2\subseteq \lambda^+$ and
$h$ a lawful mapping from $u_1$ onto $u_2$
\nl
\phantom{(g)($\alpha$)}
{\it then} for some
$M'\in K_\ap$,
$\vert M'\vert=u_2$ and $h$ is a lawful $K_\ap$-isomorphism
from
\nl
\phantom{(g)($\alpha$)}
$M$ onto
$M'$.
\nl
\phantom{(g)}($\beta$) lawful $K_\ap$-isomorphisms
preserve 
$\le_{K_\ap}$.
\item{(h)} If $ M \in K_\ap $ and $ \beta < \l^+ $ {\it then} for
some $M' \in K_\ap$ we have  $M \le_{K_\ap} M'$
and $\beta \in \vert M' \vert $
\item{(i)} [Amalgamation] Assume 
$M_\ell \in  K_\ap $  for $ \ell < 3 $  and $M_0 \le _{K_\ap} M_\ell$
for $\ell = 1,2$ . {\it Then} for some $M \in K_\ap$  and
lawful function $f$ we have:
$M_1 \le_{K_\ap} M$ , the domain of $f$ is $M_2$ , 
$f \uhr \vert M_0 \vert  $ is the     
 identity and $f$ is a $\le_{K_\ap}$-embedding of $M_2$ into $M$,
i.e. $f^{d,d}????(M_2)\le_{K_{\ap}} M$
\item{(j)} If $M_i\in AP$ for $i<i^*<\lambda$ is
$<_{K_{\ap}}$-increasing,
$\bigwedge_{i<i^*}\bigwedge_{\ell<2}M_i\le_{K_{\ap}}N^\ell\in K_{\ap}$
{\it then} there is $N^+$,
$N^2\le_{K_{\ap}} N^+\in K_{\ap}$ and a $\le_{K_{\ap}}$-embedding
$f$ of $N$ into $N^+$ over
$\bigcup_{i<i^*} M_i$.
\lz
{\bf 4.1A Remark:}
 1) This is similar to $\lambda^+$-uniform $\lambda$ forcing,
see [Sh107],
[ShHL 162] see also
[Sh326, AP],
[Sh405, AP].
\nl
2) From (g)$(\alpha), (\beta)$ we can deduce
\nl
$(\gamma)$ if $h$ is a lawful $K_\ap$-
isomorphism
from $M_1\in
K_\ap$
onto $M_2\in K_\ap$,
and $M_1\le_{K_\ap} M'_1$ 
 and $h$ can be extended to some lawful $h^+$ with domain $\vert 
M'_1 \vert$
 {\it then}
for some $h', M'_2$
we have $M_2\le_{K_\ap}M'_2$,
$h\subseteq h'$ and
$h'$ a lawful $K_{\ap}$-isomorphism from $M'_1$ onto $M'_2$.
\nl
3) We can use a linear order $<^*$ of $\lambda^+$ is
$<^*\uhr [\lambda\alpha, \lambda\alpha+\lambda)$
is a saturated model of ${\rm Th}(\bbq,<)$ and
demand legal (and lawful) maps to preserve it. No real change.
\lz
{\bf    
4.1B    
Definition }
We call $K_\ap$  {\it homogeneous} if in clause (g) of definition 4.1
we can replace ``lawful'' by ``legal''.
\lz

\lz
{\bf 4.2 Definition:}
1) For $K_\ap $ is a $\lambda$-approximation family, we let:
\nl
$K_\md=\{\Gamma:$
(i) $\Gamma$ is a $\le_{K_\ap}$-directed subset of $K_\ap$
\nl
\phantom{$K_\md=\{\Gamma :$}
(ii) $\Gamma$ is maximal
in the sense that : for every $\beta < \lambda^+$ for some
$ M \in \Gamma$ we have 
\phantom{$K_\md=\{\Gamma :$ (ii)}
$ \beta \in \vert M \vert$
\nl
\phantom{$K_\md=\{\Gamma:$}
(iii) if $M\in \Gamma$,
$M\le_{K_\ap}M'$,
 {\it then} for some $M''\in \Gamma$,
there is a lawful
\nl\phantom{$K_\md=\{\Gamma:$ (iii)}
 $K_\ap$-isomorphism $h$ from $M'$ onto $M''$ over $M\}$.
\nl
2) $K_\ap$ is a simple $\lambda$-{\it approximation} if:
(it is a $\l$-approximation family and )
{\it for every } $\Gamma\in K_{{\rm md}}$ and
$\{(M_i,N_i):i<\lambda^+\}$
satisfying $M_i\in \Gamma$,
$M_i\le_{K_{{\rm ap}}} N_i\in K_\ap$
{\it there is} a club $C$ of
$\lambda^+$ and pressing
down $h:
C\to \lambda^+$ such that:
\item{$(*)$} if $\delta_1<\delta_2$ are in $C \cap S^{\lambda^+}_\lambda$,
$h(\delta_1)=h(\delta_2)$ and $M_{\delta_1}\le_{K_{{\rm ap}}} M\in \Gamma$,
$M_{\delta_2}\le_{K_{{\rm ap}}} M\in \Gamma$ {\it then} we can find $N\in
K_{{\rm ap}}$,
$M\le_{K_{{\rm ap}}} N$,
and a lawful 
$\le_{{K_\ap}}$-embeddings $f_{\delta_1}, f_{\delta_2}$ of $N_{\delta_1},
N_{\delta_2}$ into $N$ over
$M_{\delta_1}, M_{\delta_2}$
respectively such that
$f_{\delta_1}\uhr (N_{\delta_1}\uhr \delta_1)=f_{\delta_2}\uhr
(N_{\delta_2}\uhr \delta_2)$.
\nl
Of course by refining $h$ we can demand on $\delta_1, \delta_2$ also
that
\par\noindent
\item{$(**)$}
$M_{\delta_1}\uhr\delta_1=M_{\delta_2}\uhr\delta_2$,
$N_{\delta_1}\uhr\delta_1=N_{\delta_2}\uhr\delta_2$,
$\vert M_{\delta_1}\vert\subseteq\delta_2$,
$(\vert N_{\delta_1}\vert\subseteq\delta_2$ and some $f$ is a lawful isomorphism
from $N_1$ onto $N_2$ mapping $M_1$ onto $M_2$.
\par\noindent
3) We define
$K^\alpha_\md$ as before but $M\in
\Gamma \Rightarrow \vert M\vert
\subseteq
\lambda\alpha$.
\nl
4) $K_{{\rm ap}}$ is $\theta$-closed if $\theta=\cf\theta<\lambda$
and: if
$\langle M_i:i<\theta\rangle $ is
$\le_{K_{{\rm ap}}}$-increasing in $K_{{\rm ap}}$ {\it then}
$\bigcup_{i<\theta}M_i\in K_{{\rm ap}}$ is an $\le_{K_{{\rm ap}}}$-
upper
bound; moreover
$ ( \forall i < \theta ) [M_i \le_{K_\ap} N$]
implies
$ \cup_{i<\theta} M_i  \le_{K_\ap} N $.
\nl
5) $K_\ap$ is ($<\lambda$)-closed if it is $\theta$-closed
for every $\theta < \lambda $
\nl
6) $K_\ap$ is {\it smooth} if

($\alpha$) it is ($<\l )$-closed;

($\beta$) all vocabularies $\tau_i$ are finitary;

($\gamma$) in clauses (c),(e)($\alpha$), and
  (e)($\gamma$) we can replace
``$\delta \in S^{\l^+}_\l$''    to `` $\delta >0 $ is divisible
\nl
\phantom{($\gamma$)} by $\l$''.
\nl
7) $K_{\ap}$ is a $\lambda$-approximation$^x$ family if from Definition
4.1
it satisfies clauses (a), (b), (c), is smooth??,
(g), (h), (i), and
\nl
(i)$'$ if $M_2\uhr \delta\le M_1\le\delta_1$,
then
$M_1, M_2$ have an upper bound
\nl
(j)$'$ if $M_i\le_{K_{\ap}}$ is
$\le_{K_{\ap}}$-increasing then $\bigcup_{i<\delta}M_i\le K_{\ap}$.
\nl
8) $K_{\ap}$ is nice if whenever $M_0\le_{K_{\ap}}M_1$,
$\delta\in S^{\lambda^+}_\lambda$,
$\vert M_1\vert\subseteq\delta$,
$M_0=M_2\uhr \delta$ and
$M_\ell\le M'\in K_{\ap}$ for $\ell<3$,
{\it then } we can find $M_3\le_{K_{\ap}}M''$ such that
$M'\le_{K_{\ap}}M''$,
$M_\ell\le_{K_{\ap}}M_3$ for $\ell<3$ and $M_3\uhr\delta=M_1$.
\nl
9) $K_{\ap}$ is weakly nice if whenever for $\ell=1,2$,
$M_0=M_\ell\uhr \delta_\ell$,
$\delta_\ell\in S^{\lambda^+}_\lambda$,
$\vert M_1\vert\subseteq \delta_2$, $\delta=\delta_1$
and $M'$ as above, we can find $M''$ as above.
\lz
{\bf 4.2A Observation:} 1) If $M, N\in K^\alpha_{{\rm md}}$,
$\alpha<\lambda^+$,
then  some lawful $f$ is an isomorphism from $M$ onto $N$.
\nl
2)???
\lz
{\bf 4.3 Lemma:} Suppose that
\item{(A)} $\lambda=\lambda^{<\lambda}$;
\item{(B)} $K_\ap$ is a
$\lambda$-approximation family;
\item{(C)} $\Gamma^*_\alpha\in K_{{\rm md}}$ for $\alpha<\alpha^*$;
\item{(D)} $\CT$ is a wide $\lambda^+$-tree,
$A_\alpha$ a $\lambda^+$-branch of $T$ for $\alpha<\alpha^*$
and for $\alpha\not=\beta(<\alpha^*)$ we have $A_\alpha\not=
A_\beta$,
and we
let $\vare(\alpha, \beta)=$ the level of the
$<_\CT$-last member of $A_\alpha\cap A_\beta$,
$\zeta(\alpha, \beta)=(\vare(\alpha, \beta)+1)\lambda$.

\noindent
{\it Then} there is a forcing notion $Q$ such that:
\item{(a)} $Q$ is $\lambda$-complete of cardinality $\vert
\alpha^*\vert^{<\lambda}$
\item{(b)} $Q$ satisfies the version of $\lambda^+$-c.c. from
[Sh 288 \S 1 ]   (for simplicity - here always for 
$ \vare = \w$ but by smoothness we actually have $\lub$).
\item{(c)} For some $Q$-names $\name h_\alpha$ and $\name
\Gamma'_\alpha$
(for $\alpha<\alpha^*$)
we have:
$\Vdash_Q$  ``for $\alpha<\alpha^*$ we have
$\name \Gamma'_\alpha\in K_{{\rm md}}\; ,
\name h_\alpha$ is lawful  , maps $\lambda^+$ onto
$\lambda^+$   , and maps $\Gamma_\alpha $ onto $\Gamma'_\alpha$ 
such that
for $\alpha<\beta<\alpha^*$,
$\name\Gamma'_\alpha\uhr\zeta(\alpha, \beta)=
\name \Gamma'_\beta\uhr \zeta(\alpha, \beta)$, so for every $M\in
{\Gamma}_\alpha$
we have
$h_\alpha \uhr (\vert M\vert)$ is 
lawful
and is an isomorphism from 
$M$ onto
some
$M'\in \name \Gamma'_\alpha$".
\lz
{\bf 4.3A Remark:}
 1) Our freedom is in permuting
$(\lambda\alpha,\,\lambda\alpha+\lambda)$;
up to such permutation
$\Gamma_\alpha\uhr(\lambda i)=\{ M\in
\Gamma'_\alpha:\vert M\vert\subseteq \lambda i\}$
is unique.
\nl
2) If we demand that
$K_{{\rm ap}}$
be smooth
the proof is somewhat simplified.
\nl
3) We can replace assumption (B) by
\nl
 (B)$'$ $K_{\ap}$ is a
$\lambda$-approximation$^x$ family.
\hz
{\it Proof:}
We define $Q$ as follows:
\nl
$p\in Q$
{\it iff} $p=\langle (M^p_\alpha, h_\alpha^p ) :\alpha\in w^p\rangle$
where

(a) $w_p\in [\alpha^*]^{<\lambda}$;

(b) $M^p_\alpha\in \Gamma_\alpha$;

(c) $h^p_\alpha$ a 
lawful
 mapping,
Dom $h^p_\alpha=\vert M^p_\alpha\vert$;

(d)  if
$\alpha\not=\beta$ are in $w^p$,
then  :
$h_\alpha(M^p_\alpha\uhr \zeta(\alpha, \beta))$
and
$h_\beta(M^p_\beta\uhr \zeta(\alpha, \beta))$ are
$\le_{K_\ap}$-
\nl
\phantom{(d) (d)}
comparable;

(e) for every $\alpha\in w^p$,
for some $n<\omega$,
$0=i_0<i_1<\dots < i_n=\l^+$,
we have: for
\nl\phantom{(e) (e)}
 $\ell \in [1,n)$,
$i_\ell \in S^{\l^+}_\l$
and for every $\ell <n$
\nl
\phantom{(e) (e)}
$(*)_\ell $ for every
$\beta \in w$ for which
$ \zeta (\alpha, \beta ) \in [i_\ell , i_{\ell +1 } )$
 and $ j \in   [i_\ell , 
i_{\ell + 1 } ) 
\cap S^{\l^+}_\l$
there is
\nl
\phantom{(e) (e) $(*)_\ell $}
 $ \gamma \in w$  such that:
$j \le \zeta 
( \alpha , \gamma )
  \in [i_\ell , i_{\ell +1 }  ) $
and
$ M^p_\beta \uhr \zeta ( \alpha , \beta ) \le_{K_\ap}
M^p_\gamma \uhr \zeta ( \alpha , \gamma ) $
\hz
{\it The order is }
$p\le q$ iff:
$w^p\subseteq w^q$ and for $\alpha\in w^p$ :
$M^p_\alpha\le_{K_\ap}M^q_\alpha$, 
$h^p_\alpha\subseteq h^q_\alpha$ and
$M^p_\alpha\not=M^q_\alpha\Rightarrow \bigwedge_{\beta\in w^p}
h_{\beta}(M^p_\beta\uhr\zeta(\alpha, \beta))\le_{K_\ap} 
h_{\alpha}(M^q_\alpha\uhr
\zeta(\alpha, \beta))$.
\nl
The lemma will follow from the facts 4.4-4.7 below.
\hz
{\bf 4.4 Fact:} Any increasing chain in $Q$ of length $<\lambda$ has an
upper bound.
\hz
{\it Proof:} Let
$\langle p_i:i<\delta\rangle$ be an increasing sequence in $Q$,
$\delta<\lambda$ a limit ordinal.
Let $w=\bigcup\{w^{p_i}:i<\delta\}$,
and list $w$ as $\{\alpha_j: j< j^*\}$.
We now choose by induction on $j<j^*$, a member $M_{j}$ of
$K_{{\rm
ap}}$ and a lawful mapping $h_j$ with domain
$\vert M_j \vert$ such that :
\item{$\otimes$}(a)  if $\langle (M^{p_i}_{\alpha_i}, h^{p_i}_{\alpha_i}):
i<\delta$ but $\alpha_j\in w^{p_i}\rangle$ is eventually constant, then
this value is
\nl
\phantom{a)} $(M_j, h_j)$.
\nl
(b)  Otherwise  let
$h_j(M_j)\in \Gamma_{\alpha_j}$ be
a
$\le_{K_{{\rm ap}}}$-upper bound of
$\{h_{\alpha_i}^{p_i}(M^{p_i}_{\alpha_i} ):\; i<\delta$ but
\nl
\phantom{(b)} $\alpha_j\in w^{p_i}\}\cup  \{
h_{j_1}(M_{j_1})\uhr\zeta(\alpha_j, \alpha_{j_1}): j_1<j\}$.
\hz

If we succeed $q=^{df}
\langle ( M_j, h_j ):j\in w\rangle$ is a member of
$Q$  as 
required.
Why?  First  we check that  $q \in Q$ . Clauses (a),(b),(c) are obvious;
for clause (d) let
$\alpha \not= \beta$ be in $w$ , so let
$\{ \alpha , \beta \} = \{ \alpha_{j_1}, \alpha_{j_2} \}$,
$j_1 < j_2 $; now if ($*$)(b) holds for $j_2$ just note that
$h_{j_1} (M_{j_1} ) \uhr \zeta ( \alpha_{j_1} , \alpha_{j_2} ) \le
h_{j_2} (M_{j_2} )$
by the choice of the later; and if ($*$)(a) holds for $j_2$,
 then for some
$i<\d , ( M_{j_2} , h_{j_2} ) = (M^{p_i}_{j_2} , h^{p_i}_{j_2} )$
and now check the choice of
$ (M_{j_1} , h_{j_1} ) $.
%
 If for it too clause (b) holds for some $i(1)<\delta$,
$(M_{j_1}, h_{j_1})=
(M^{P_{i(1)}}_{j_1}$,
$h^{P_i(1)}_{j_2})$ and use
$p_{\max\{i(1), i\}}\in Q$.
If for $j_1$ clause (b) holds then by its choice
$h^{p_i}_{j_2}(M^{P_i}_{j_2})\uhr \zeta(j_1, j_2)\le h_{j_2}(M_{j_1})$
hence $h_{j_1}(M_{j_2})\uhr\zeta(j_1, j_2)=
h^{p_i}_{j_1}(M^{p_i}_{j_1})\uhr\zeta(j_1, j_2)\le
h_{j_1}(M_{j_2})\uhr\zeta(j_1, j_2)$ as required.
So we are left with the case clause
(b) of ($*)$ apply to $j_2$, which is even easier.
For clause (e) ,  clearly it is enough to prove :
\nl
\item{$(*)$}  for every $i_1
 \in ( S^{  \lambda^+  }_ \lambda
\cup  \{        
 \lambda^+ \} ) $  
there is 
$i_0 \in i_1  
\cap  
( S^{\l^+}_\l \cup \{  
0
\} )$
such that $(*)_\ell$ of clause (e) of the definition of 
$Q$ holds with $ i_0 , i_1 $ taking the role of 
$i_\ell , i_{\ell +1 }$ .
\lz
Let $ i_1  \in S^{\l^+}_\l
\cup  \{ \l^+ \} $   be given ; for each $ i  < i_1 $  let
$f(i) =^{df}  
\sup \{ \zeta ( \beta , \alpha ) +1: \beta \in w , \zeta(
\beta , \alpha    )  
 \in [ i , i_1 ) \}$
( if the supremum is on an empty set - we are in a trivial case ).
Clearly 
$ [ j_1  < j_2  < i_1 \Rightarrow f(j_1 ) \geq f( j_2 ) ]$,
so for some $i_0 \in  i_1 \cap ( S^{\l^+}_\l  \cup  \{ 0 \} )$
for all  $i \in [i_0 , i_1 ) \cap  ( S^{\l^+}_\l \cup \{0\})$
we have  $ f(i ) = f( i_0 )$ .
Now 
for each $ i < i_1 $
let
$g(i) =^{df}  \sup \{  j +1: j < j^* , \zeta ( \alpha_j , \alpha )
\in  [i , i_1 )$ and in $\otimes$ case (b) occurs for $j\} $,
note: 
 if the suremum is on the empty set then the value
is zero;
again it is clear that $g$ decrease with $i$
hence wlog for
all  $ i \in [i_0 , i_1 ) $ we have
$g(i) = g( i_0 ) $
\nl
{\bf case 1} $g(i_0)>0$; this means that for every $i \in [i_0 , i_1 )$ there is
$\beta \in w$ such that  :
$\zeta ( \beta , \alpha ) \in [ i_0 , i_1 )$
and  
letting
 $\beta = \alpha_j$ and in (*) above case (b) occurs.
 \nl
Check
\nl
{\bf case 2} not case 1
\nl
For every $ \gamma \in w$  let  
$\xi_\gamma$ be the
first ordinal $ \d $ such that 
$\lk ( M^{p_i}_\gamma  , h^{p_i}_\gamma ) : i < \d, i\geq \xi_\gamma \rk$
is constant, and again wlog  for some 
$\varepsilon^*$
for 
every
$ i \in [ i_0 , i_1 ) , \varepsilon'< \varepsilon^*,
\zeta < f( i_0 )$
and $j < g( i_0 ) $ there is $\beta \in w$
such that   
 $ \zeta \le \zeta ( \beta , \alpha ) \in [ i, i_1) , 
\beta \in \{ \alpha_{j' }:
j \le j' < g_a( i_0 ) \}$
and $j_\beta \ge \varepsilon' $,
the
 rest should be clear.

So we have proved that $q \in Q$ ; now $p_i \le_{K_\ap} q$
is straightforward. So now we have only to prove that we can carry 
the inductive definition from (*).

\hz
In the choice of $M_j,$ $h_j$ we first have chosen $h_j(M_j)$. We do it 
by
choosing $h(M_j\uhr \zeta)$ for $\zeta\in\{\zeta(\alpha_j,
\beta):\beta\in
w\}$;
there 
we use clause (e)$(\d   )$ of Definition 4.1.
Having chosen $h_j(M_j)$ we can find $M_j$, $h_j$ by clauses
(g)($\alpha)+(\beta)$ of Definition 4.1.
\hfill\qed$_{4.4}$
\lz
{\bf 4.5 Fact:} 1) if $p\in Q$,
$\alpha\in w^p$  and
$N\in \Gamma_\alpha$ 
then for some $q$:  
$p\le q$,
$w^q=w^p$  and
\nl
$\bigwedge_{\beta\in w^p\setminus\{\alpha\}}
(M^p_\beta, h^p_\beta)=(M^q_\beta, h^q_\beta)$ and
$N\le M^q_\alpha$.
\nl
2) If $p\in Q$,
$\alpha<\alpha^*$ then for some $q$,
$p\le q\in Q$ and $\alpha\in w^q$.
\hz
{\it Proof:} 1) Easier than the previous one
(or  
 let $\d =1 $ , $p_0 = p$ and 
$\{ \alpha_j : j < j^* \} $ list  $w^p$
with $ \alpha = \alpha_0$, repeat the proof of 
4.4 , just use $q$ to choose $ ( M_0 , h_0 ) $.\hfill\qed$_{4.5}$
\nl
2) Easier.
\hz
Note the following	
\lz
{\bf 4.6 Fact:} If $K_{{\rm ap}}$ is $\theta$-closed,
then the following set is $Q'$  dense in $Q$:
 $\{p\in Q$: if $\alpha, \beta\in w^p$, then
$h^p_\alpha(M^p_\alpha)\uhr\zeta(\alpha, \beta)=
h^\ell_\alpha(M^p_\beta)\uhr\zeta(\alpha, \beta)\}$.
\hz
{\it Proof:} Follows easily from the previous Facts.
\lz
{\bf 4.7 Fact:} The chain condition $(*)_{\lambda^+, \omega}$
from [Sh 288 \S1] holds.
\hz
{\it Proof:} For simplicity assume $K_{\ap}$ is $\aleph_0$-closed so we
can use 4.6. Suppose $p(\delta)\in Q$ for $\delta\in
S^{\lambda^+}_\lambda$.
For some pressing down function $h: S^{\lambda^+}_\lambda\to 
\lambda^+$ and $\langle \omega_\gamma:\gamma<\lambda^+\rangle$ we
have:
\item{$(*)$}
if $h(\delta^1)=h(\delta^2)$,
$\delta^1<\delta^2$ {\it then}:
\nl
(a) $\otp (w^{p(\delta^1)})=\otp (w^{p(\delta^2)}$) and
$w^{p(\delta^1)}\cap w^{p(\delta^2)}=w_{h(\delta^1)}$
\nl
(b${}_1$) $OP_{w^{p(\delta^1)}, w^{p(\delta^2)}}$
is\footnote{$^\dagger$}{$OP_{u_1, u_2}$ is the unique order preserving
function $f$ such that  $(\Dom f)$ an initial segment of $u_1$ $\Rang(f)$ an initial
segment of $u_2$ and $\Dom(f)=u_1\vee\Rang(f)=u_2$.} the identity on
$w^{p(\delta^1)}\cap w^{p(\delta^2)}$ 
\nl
(b${}_2$) for $\alpha,\beta\in w^{p(\d^1)}$ the following are 
equivalent:
\nl
\phantom{(bbb)}(i) $\zeta(\alpha,\beta)<\delta^1$;
\nl
\phantom{(bbb)}(ii) $\zeta(\alpha',\beta')<\d^2$ where $\alpha'=^{df}
\OP_{w^{p(\d^1)},w^{p(\d^2)}}(\alpha)$, $\beta'=^{df}
\OP_{w^{p(\d^1)},w^{p(\d^2)}}(\beta)$;
\nl
\phantom{(bbb)}(iii)
$\zeta(\alpha',\beta')=\zeta(\alpha,\beta)$ where $\alpha'$,
$\beta'$ are as in (ii).
\nl
(c) $\vert M^{p(\delta^1)}_\alpha\vert$ is bounded in $\delta_2$
and also
$ \sup \{ \zeta ( \alpha , \beta )$ : $\alpha \not= \beta $ are
in   $w^{p(\delta^1)} \}
< \d_2$
\nl
(d) if $\alpha^2=\OP_{w^{p(\delta^1)}, w^{p(\delta^2)}}(\alpha^1)$ then
\nl
\phantom{(d)} ($\alpha)$ $\OP_{\vert
M^{p(\delta^1)}_{\alpha^1}\vert,
\vert M^{p(\delta^2)}_{\alpha^2}\vert}$ is an isomorphism from
$M^{p(\delta^1)}_{\alpha^1}$ onto
$M^{p(\delta^2)}_{\alpha^2}$ which is lawful.
\nl
\phantom{(d)}
($\beta$) $M^{p(\delta^1)}_{\alpha^1}\uhr (\delta^1\lambda)=
M^{p(\delta^2)}_{\alpha^2}\uhr( \delta^2\lambda)$.
\lz
Now we  have
to prove  :
$ h( \delta^1 )= h( \d^2 ) \Rightarrow   p(\d^1) , p( \d^2 ) $
are compatible. In the list
$\{ \alpha_j : j < j^* \} $
 put
$w^{p(\d^1 )} \cap  w^{p( \d^2 ) } $ an initial segment.
 Say $\{\alpha_j:j<i^*\}$.
First
we  restrict ourselves further by assuming
$K_{\ap}$ is nice (see Definition 4.2(7)).
We define a common upper bound $p$;
we let
$w^p=w^{p(\delta^1)}\cup
w^{p(\delta^2)}$.
For $\alpha\in w^{p(\delta^\ell)}\setminus w^{p(\delta^{1-\ell})}$ let
$(M^p_\alpha, h^p_\alpha)=
(M^{p(\delta^\ell)}_\alpha,
h^{p(\delta^\ell)}_\alpha)$. For $\alpha\in w^{p(\delta^\ell)}\cap
\omega^{p(\delta^\ell)}$
first choose
$M^p_\alpha\in \Gamma_\alpha$ such that
$M^{p(\delta^\ell)}_\alpha\le_{K{\ap}}M^p_\alpha$ for
$\ell=1, 2, $ and
$M^p_\alpha\uhr\delta^1=M^{p(\delta^\ell)}_\alpha\uhr\delta^2$
[Why? -by Definition 4.1 clause (e)($\gamma$), now we can find such
$M^p_\alpha\in K_{\ap}$, now we can find one in $\Gamma_\alpha$ by
$``K_{\ap}$ is nice'' (see Definition 4.2(7))].
\nl
Second, we deal with the case $K_{\ap}$ is not nice. Without loss of
generality there is $\delta_0\in S^{\lambda+}_\lambda$,
$M^{p(\delta^\ell}_{\alpha_j}\uhr \delta_1\subseteq \delta_0$ and
$A_{\alpha_j}\uhr \delta_0\neq A_{\alpha_i}\uhr \delta_0$ for
$j<i<i^*$.
We choose by induction on $j\le i^*$ a condition $q^j$, increasing with
$j$,
$w^{q^j}=\{\alpha_i:i<j\}$,
$q^j$,
$\bigwedge_{\ell=1, 2}\bigwedge_{i<j}
M^{p(\delta^\ell)}_{\alpha_i}\le_{K_{\ap}} M^{q_j}_{\alpha_i}$.
The bookmkeeping is as in the proof of 4.4,
the successor case as in the proof above (for nice $K_{\ap}$) but using
amalgamation (=clause (i) of Definition 4.1 in the end).
\hfill\qed$_{4.7}$

This
 finishes the proof of 4.3.
\hfill\qed$_{4.3}$

 The simplicity of $K_{\ap}$ is referred to only in 4.8 below,
but it is needed to get the universality results later.
\lz
\lz
{\bf 4.8 Claim:}
Assume $K_{{\rm ap}}$ is a simple $\lambda^+$-approximation
system.
If $\Gamma_0\subseteq K_{{\rm ap}}$ is directed and
$\alpha<\lambda^+\Rightarrow
\lambda=
\vert [\lambda\alpha, \lambda\alpha+\lambda)\setminus \bigcup_{M\in
P}\vert M\vert$, {\it then} for some
forcing
notion
$Q$ satisfying the $\lambda^+$-c.c. of [Sh 288 \S1],
$\vert Q\vert=\lambda^+$,
$\Vdash_Q$ ``there is a $\Gamma$ and a lawful $f$ such that
$f(\Gamma_0)\subseteq \Gamma\in K_{{\rm md}}$''.
\hz
{\it Proof:}
Natural. By renaming,
 without loss of generality $A\eqdf \cup\{\vert M\vert: M\in
\Gamma_0\}=
\{2\alpha:\alpha<\lambda^+\}$.
$Q=\{M: M\in K_{{\rm ap}}$ and $M\uhr A\in K_{{\rm ap}}$ and
$M\uhr A\le_{K_{{\rm ap}}}M\}$ order by
$\le_{K_{{\rm ap}}}$.
\hfill\qed$_{4.8}$
\lz
{\bf 4.9 Conclusion:} Assume
$\lambda=\lambda^{<\lambda}<2^{\lambda^+}=\chi$ ,
and  a  $\lambda^+$-tree $\CT$ with $\ge \chi$
branches is given\footnote{$^{\dagger\dagger}$}{if $\lambda=\lambda^{<\lambda},
2^\lambda=\lambda^+<\chi=\chi^{\lambda^+}$,
and we add $\chi$ Cohen subsets to $\lambda^+$ (i.e. force by $\{f:f$ a
partial function from $\chi$ to $\{0, 1\}$ of cardinality
$<\lambda^+\}$, then in $V^P$ those assumptions hold.}
For simplicity we assume that $\lambda^+$ is the set of members of 
$\CT$,
$0$ is the root and $\alpha<_\CT\beta\Rightarrow\alpha<\beta$
for $t\in \CT$ and 
 let $u_t=\{[\alpha\lambda,
\alpha\lambda+\lambda):\alpha\le_T t\}$ .
 {\it Then} there is a forcing notion $P$ such that:
\item{(a)} $P$ is $\lambda$-complete, satisfies the $\lambda^+$-c.c.
and has cardinality $\chi$
(so the cardinals in $V^P$ are the same and cardinal arithmetic
should be clear).
\item{(b)} for any  $\l$-approximation
system $K_\ap $
there are $\langle \Gamma^\zeta_t,\, M_t: t\in
\CT_\gamma\rangle$
for
$\zeta<\lambda^{++}$ such that:
\itemitem{($\alpha)$} $\Gamma^\zeta_t
        \in K_\md^{\lambda(\lg(t)+1)}$
\itemitem{($\beta$)} $t<_Ts\Rightarrow
\Gamma^\zeta_t\subseteq
\Gamma^\zeta_s$
\itemitem{($\gamma$)} for every $\Gamma\in K_\md$ for some
$\zeta<
\l^{++}$ and $\lambda^+$-branch
$B=\{t_\alpha:\alpha<\lambda^+\}$
 of $T$
and 
lawful function from $\lambda^+$
onto $\lambda^+$
mapping $\Gamma$ onto
$\bigcup_{\alpha<\lambda^+} \Gamma_{t_\alpha}$.
\item{(c)} Is $R\in V^P$ is
$(<\lambda)$-complete, satisfies the version of the $\lambda^+$-c.c.$(*)_{\lambda^+ ,
\omega}$
 from [Sh288\S1] and
 $D_i\subseteq(i<\lambda^+)$ is a dense subset of $R$ and $\vert
R\vert\le\lambda^{+}$,
 {\it then} for some directed $G\subseteq R$,
$\bigwedge_i D_i\cap G\not=\emptyset$.
\hz
{\it Proof:}
We use iterated forcing of length $\chi \times\lambda^{++}$,
$(<\lambda)$-support , 
each iterand
 satisfying the $\lambda^+$-c.c.$(*)_{\lambda^+,\omega}$ from [Sh
288 \S1],
$\langle P_i, \name Q_j:
i\le \chi+\lambda^{++}, j<\chi \times \lambda^{++}\rangle$
 such that: for every $K_\ap$
(from $V$ or
 from some intermediate universe)
for unboundedly many $i<\chi\times \lambda^{++}$,
we use the forcings
 from 4.3
or 4.8.
\hfill\qed$_{4.9}$
\bigskip
\noindent
{\bf \S 5 Applications}
\medskip\noindent
{\bf 5.1 Lemma:}
Suppose
\item{(A)} $T$ is first order , complete, for simplicity with elimination
of quantifiers (or just
inductive theory with the amalgamation and
disjoint embedding property).
\item{(B)} $K_\ap $ is a simple $\lambda$-approximation system such
that
every $M\in K_\ap$ is a model of $T$ hence every $M_\Gamma$ ,
where for
$\Gamma\in K_\md$ we  let
$M_\Gamma=\bigcup\{M:M\in \Gamma \}$.
\item{(C)} every model $M$ of $T$ of cardinality $\lambda^+$
can be embedded into $M_\Gamma$
\nl
for some $\Gamma \in K_{\md}$ with $\bigcup_{M\in P}\vert
M\vert=\{2\alpha:\alpha<\lambda^+\}$.
\nl

\noindent
{\it Then:}

\item{(a)} in 4.9 in $V^P$,
there is a model of $T$ of cardinality $\lambda^{++}$ universal for
models of $T$ of cardinality $\lambda^+$.
\nl
\item{(b)} So in $V^P$ , $\univ (\l^+ , T ) \le \l^{++}$ but there is a
club guessing sequence
$\langle C_\delta:\delta\in S^{\lambda^+}_x\rangle$.
\lz
{\it Proof:} Straightforward.
\hfill\qed$_{5.1}$
\lz

Though for theories with the strict order property,
the conclusion of \S4 (and 5.1) fails,
for some non simple theories we can succeed.
Note that in
5.1  we have some freedom in choosing $K_\ap$ even after $T$ is
fixed .
\lz
{\bf 5.2 Lemma:}
Let $T=T^*_{\feq}$; it satisfies the assumption of 5.1 (hence its
conclusions).
\nl
In fact we can find a smooth nice
simple
$\l$-approximation system
$K_\ap$  such that every model $M$ of $T$ of cardinality $\l^+$ is
embeddable into some
$M \in K^\md_{\ap} $ .
\hz
{\bf 5.2A Remark}
1) Note that there  $\univ (\l , T^*_\feq ) = \univ ( \l , T_\feq )$ .
Actually the $\l$-approximation family we get is also homogeneous.
\nl
2) The situation is similar for $T_3$  in 5.3.
\hz

{\it Proof:}
By 5.2A(1) we deal with models of $T_{\feq}$.
Condition (A) of 5.1 clearly holds.
The main point is to define $K_\ap$.
\item{($\alpha$)} $M\in K_\ap$ iff:
\itemitem{(i)} $M$ is a model of $T$
\itemitem{(ii)} $\vert M\vert\in [\lambda^+]^{<\lambda}$
\item{($\beta$)} $M_1\le_{K_\ap} M_2$ iff
\itemitem{(i)} $M_1\subseteq M_2$
\itemitem{(ii)} if $\delta\in S^{\lambda^+}_\lambda$,
$a\in P^{M_1}\cap\delta$,
$b\in Q^{M_1}\setminus \delta$ and
$(\forall c\in M_1)
[M_1\models ``bE_a c$''$\Rightarrow c\notin \delta]$
{\it then }
$(\forall c\in M_2)$
$[M_2\models ``b E_ac$''$\Rightarrow c\notin \delta]$.
\lz
Also condition (C) of 5.1 is easy and we turn to condition (B).
 The checking of ``$( K_\ap, \le_{K_\ap} ) $ is a
$\l$-approximation family'' (see Definition 4.10)
as well as  smoothness
 is staightforward.
E.g. let us check the amalgamation (Definition 4.1 clause(i)).
 So assume $M_\ell\in K_{\ap}$ for $\ell<3$,
$M_0\le_{K_{\ap}}M_1$,
$M_0\le_{K_{\ap}}M_2$;
by Definition 4.2  clause (g)$(\alpha)$ without loss of generality
$\vert M_1\vert\cap \vert M_2\vert=\vert M_0\vert$. Now we shall
 define a model $M$
with universe  $\vert M_1 \vert \cup  \vert M_2 
 \vert$,
as follows:
$ P^M  =^{df} P^{M_1} \cup  P^{M_2}$ , 
$ Q^M  =^{df} Q^{M_1} \cup  Q^{M_2}$ , 
and for each $x \in P^M$, we let
$E_x$ be the
closure to an equivalence relation of
the
 set of cases occurring in $M_1$ and/or $M_2$, now check.
The checking is straightforward.
\hz
Now we are left with the main point: the simplicity of $K_{\ap}$
 (see Definition 4.2(2)). Choose $h$ as implicit in
$(**)$ of Definition 4.2(2);
so let $\delta_1<\delta_2$,
$M_{\delta_1}, M_{\delta_2}, M, N_{\delta_1}, N_{\delta_2}, f$ be as there.
Let
$f_{\delta_1}, f_{\delta_2}$ be  lawful mappings such that
$f_{\delta_1}\uhr(N_{\delta_1}\uhr\delta_1)=
f_{\delta_2}\uhr(N_{\delta_2}\uhr\delta_2)$
and $N'_0\eqdf M$
and for $\ell=1, 2$
 the
 $f_{\delta_\ell}(N_{\delta_\ell}\setminus M_{\delta_\ell})$ is disjoint to $M$
for $\ell=1, 2$;
let
$N'_\ell=f_{\delta_\ell}(N_{\delta_\ell})$.
Now we define $N\in K_{\ap}$; it is a model with universe $\vert
N'_0\vert\cup
\vert N'_1\vert\cup\vert  N'_2\vert$,
$P^N=p^{N'_0}\cup P^{N'_1}\cup P^{N'_2}$,
$Q^N=Q^{N'_0}\cup Q^{N'_1}\cup Q^{N'_2}$. Lastly
for $x\in P^N$,
we let $E_x$ be the finest equivalence relation on $Q^N$ which extend
each
$E^{N'_\ell}_x$,
(if $x\in N'_\ell$).
Why is $N$ a model of $T_{\feq}$?
Clearly $P^N, Q^N$ is a partition of $N$ (as this holds for
$N'_\ell(\ell<3)$ and as any two of those models agree on their
intersection) and each
$E^N_x(x\in P^N$) is an equivalence relation on $Q^N$ (by its choice).
Why $N'_\ell$ is a submodel of $N$?
concerning $P$ and $Q$ there are no problems.
So assume $x\in P^{N'_\ell}$,
and we shall prove
$E^N_x\uhr\vert N'_\ell\vert=E^{N'_\ell}_x$,
the inclusion $\supseteq$ is by the choice of $E^N_x$.
For the other inclusion, if $x$ belongs (and the proof of amalgamation) to only one $N'_m$ it is
totally trivial.
If it belongs to exactly two of them,
say $N'_{m_1}$, $N'_{m_2}$ just note $
N'_{m_1}\cap N'_{m_2}\le N'_{m_1}, N'_{m_2}$.
So assume $x\in \bigcap^2_{\ell=0}P^{N'_\ell}$ and here we shall use
clause $(\beta)$(ii) of the definition of $\le_{K_{\ap}}$.
So suppose
$y_0, \dots, y_{m(*)}$ are such that
$N'_{\ell(*)}\models y_mE_xy_{m+1}$,
but
$N'_{\ell(m(*))}\models \neg y_0E_xy_{m(*)}$; wlog $m(*)$ is
minimal.
Of course without loss of generality $\ell(m)\neq \ell(m+1)$
(as then we can omit $y_{m+1})$
and $\ell(m)\neq \ell (m+2)$
(otherwise $y_{m+1}\in N'_{\ell(m)}\cap
N'_{\ell(m+1)}=N'_{\ell(m+1)}\cap N'_{\ell(m+2)}\to y_{m+2}$ and we can
omit $y_{m+1}$.
So necessarily
$\ell(m)=\ell(m+3)$ and $\{\ell(m), \ell(m+1), \ell(m+2)\}=\{0, 1, 2\}$,
hence enough to deal with the case $m(*)=3$.
As $x\bigcap_{m<\zeta}N'_{\ell(m)}=??? \in \bigcap_{\ell=0}^2N'_\ell$,
clearly $x\in N'_1\cap N'_2\subseteq \delta_1$, and for some $m$,
$y_m\in N'_1\cap N'_2(\subseteq \delta_1)$, and for some $m_1\in \{m-1,
m+1\}$ and
$k\in \{1, 2\}$,
$y_{m_1}\in N'_0\cap N'_k=M\cap N'_k=M\cap N_k$,
so by the choice of $\le_{K_{\ap}}$ there is $y'_{m_1}\in N'_0\cap
N'_1\cap \delta_k\cap Q^N$,
$y'_{m_1}E_x^{N_0'}y_m$, but so $y'_{m_1}\in \bigcap_{\ell<3}N'_\ell$
(as $N'_1\cap \delta_1=N'_2\cap \delta_2$), and we are done.
\hfill\qed$_{5.2}$
\lz
{\bf 5.3 Lemma:}
$T_{{\rm trf}}$,
the theory of triangle free graphs satisfies the assumption of 5.1
(hence its conclusions).
\hz
{\it Proof:}
Let $xRy$ mean $\{x, y\}$ is an edge.
 The main point is to define $K_\ap$
\item{($\alpha$)} $M\in K_\ap$ iff
\itemitem{(i)} $M$ is a model of $T$
\itemitem{(ii)} $\vert M\vert \in [\lambda^+]^{<\lambda}$
\item{$(\beta)$} $M_1\le_{K_\ap}M_2$ iff
\itemitem{(i)} $M_1\subseteq M_2$
\itemitem{(ii)} if $\delta\in S^{\lambda^+}_\lambda$,
$a,b \in M_1\cap \delta$ and there is no
$c\in M_1\cap \delta$,
$M_1\models cRa\,\&\,cRb$
{\it then} for no $c\in M_2\cap \delta$,
$M_2\models cRa\,\&\, cRb$.

Let us check Definition 4.1, i.e.
that
$(K_{\ap}$, $\le_{K_{\ap}})$ is a $\lambda$-approximation system
\nl
{\it Clause (a), (b), (c)} are immediate.
\nl
{\it Clause (d)} holds in a strong form: the natural union is a $\lub$;
and even $K_{\ap}$ is smooth.
\nl
{\it Clause (e)} follows from (d)$^+$ and (i)
(amalgamation)
\nl
{\it Clauses (f) and (g)} are immediate (as in g) the demand is on
lawful $h$ only).
\nl
{\it Clause (h)} is trivial.
\nl
{\it Clause (i):} Using a lawful $f$ without loss of generality $\vert
M_1\vert\cap
\vert M_2\vert=\vert M_0\vert$.
Define $M_3:\vert M_3\vert=\vert M_1\vert \cup\vert M_2\vert$,
$R^{M_3}=R^{M_1}\cup R^{M_2}$.
Clearly $M_3\in K_{\ap}$ as for $M_\ell\le_{K_{\ap}}M_3$, by transitivity
and symmetry it is enough to prove $M_1\le_{K_{\ap}}M_3$, clearly
$M_1\subseteq M_3$,
(i.e. clause (i) of ($\beta$) above).
For proving clause (ii) let
$\delta\in S^{\lambda^+}_\lambda$, $a, b\in M_1$ and
$c\in M_3\cap \delta$
be such that $aR^{M_3}c\&bR^{M_3}c$.
If $c\in M_2\setminus M_1$, necessarily $a, b\in M_2$ hence $a,b\in M_0$
and use $M_0\le_{K_{\ap}} M_2$, but if $c\in M_1$ there is nothing to
prove.
\nl
{\it Clause y} follows from smoothness.

Next let us show that $K_{\ap}$ is simple.
Let $\delta_1<\delta_2$ (from $S^{\lambda^+}_\lambda)$,
$M_{\delta_1}, M_{\delta_2}, N_{\delta_1}, N_{\delta_2}, M,f$ be as in
Definition 4.2(2) $(**)$.
Without loss of generality $M\cap N_{\delta_1}=M_{\delta_1}$, $M\cap
N_{\delta_2}=M_{\delta_2}$. Define a model $N$.
$$\vert N\vert=\vert N_{\delta_1}\vert \cup\vert N_{\delta_2}\vert \cup
\vert M\vert$$

$$R^N=R^{N_{\delta_1}}\cup R^{N_{\delta_2}}\cup R^{N_{\delta_3}}$$

Clearly $N$ extends each of the models $N_{\delta_1}, N_{\delta_2}, M$
 (hence $M_{\delta_1}, M_{\delta_2}$ too).

Clearly for proving $N\in K_{\ap}$ it suffices to show
\item{$(*)_1$} if $N$ there is no triangle.

\noindent
Why?
Clearly the only case we should consider is $a\in N_{\delta_1}\cap
N_{\delta_2}\setminus M$,
$b\in N_{\delta_1}\cap M\setminus N_{\delta_2}$,
$c\in N_{\delta_2}\cap M\setminus N_{\delta_1}$.
(hence $b\in M_{\delta_1}\setminus \delta_1$,
$c\in M_{\delta_2}\setminus \delta_2$).
So for some $c'\in M_{\delta_1}\setminus \delta_1$,
$f(c')=c$ but also $f(a)=a$ hence $aR^{N_{\delta_1}}c^1$, also
$aR^{M_{\delta_1}} b$, so as
$M_{\delta_1}\le_{K_{\ap}}N_{\delta_1}$ clearly for some
$a'\in M_{\delta_1}\cap \delta_1$ we have
$a'R^{M_{\delta_1}}c' R a' R^{M_{\delta_1}}b$. Applying again $f$ we get
$a'R^{M_{\delta_2}}c$. So (by the last two sentences) in $M$ we have
$a'R^M c\&a'R^M b$. But by the choice of abc (and as $M\subseteq N$) we
get $a',b,c$ is a triangle in $M$ which belongs to $K_{\ap}$,
contradiction.

So $N\in K_{\ap}$; also $M\le_{K_{\ap}}N$.
[Why? being submodels sould be clear. So suppose $\delta, a, b, c$
contradicts clause (ii) of the definition of $\le_{K_{\ap}}$, so $c\in
N\setminus M$,
so $c\in N_{\delta_\ell}\setminus M_{\delta_\ell}$ for some $\ell\in
\{1, 2\}$.

If $c\in N_{\delta_\ell}\setminus M_{\delta_2}\setminus
N_{\delta_{1-\ell}}$ then necessarily (as $aR^Nc, bR^Nc)$ we have $a, b\in
N_{\delta_\ell}$ hence $a,b\in N_{\delta_\ell}\cap M=M_{\delta_\ell}$ using
$M_{\delta_\ell}\le_{K_{\ap}} N_{\delta_\ell}$ we have
$c'\in M_{\delta_\ell}\cap \delta$ such that $M_{\delta_\ell}\models
aRc'\&bRc'$,
 $c'$ is as required).

So necessarily $c\in N_{\delta_\ell}\cap N_{\delta_{1-\ell}}\setminus M$
hence $a, b\in M_{\delta_1}\cup M_{\delta_2}$;
if $a,b\in M_{\delta_1}$ do as above, also if $a,b,\in M_{\delta_2}$ do
as above, so by symmetry without loss of generality $a\in M_{\delta_1}$,
$b\in M_{\delta_2}$. Now use $f$ as in the proof of ``$N$ is triangle
free''.)

The case $N_{\delta_\ell}\le N$ is similar.

Having proved $K_{\ap}$ is a simple $\lambda$-approximation family, it is
easy to check the assumptions (A) and (C) of 5.1 hold.
\hfill\qed$_{5.3}$
\lz
{\bf 5.4 Discussion:}
The similarity between the proofs of 5.2, 5.3 is not incidental.
For a complete first order $T$,
let
e.g. $M^*\in \GB(\chi)$ be a $\lambda^+$-saturated model of $T$,
choose by induction on $\zeta<\lambda$, an elementary submodel
$\GB_\zeta$ of
$(\GB(\chi), \in , <^*_\chi$) of cardinality $\lambda^+$ such that
$\{M^*, \GB_\varepsilon:\varepsilon<\zeta\}\cup(\lambda^++1)\subseteq
\GB_\zeta$,
$\GB_\zeta^{<\lambda}\subseteq \GB_\zeta$,
$f_\zeta\in \GB_{\zeta+1}$ a mapping from $\GB_\zeta\cap \lambda^+$ onto
$\{\delta+\vare:\vare\le\zeta\delta<\lambda^+$ divisible by $\lambda\}$,
extending $\bigcup_{\varepsilon<\zeta}f_\varepsilon$.
In the end let $N^*$ be the model with universe $\lambda^+$ such that
$\bigcup_{\zeta<\lambda^+}f_\zeta$ as an isomorphism for
$\bigcup_{\zeta<\lambda^+}\GB_\zeta\cap M^*$ onto $N^*$.
Let
$E$ be a thin enough club of $\lambda^+$.
Let
$K'_{\ap}=\{N:N\prec N^*$ and
$(N, \delta)_{\delta\in N\cap \nacc E}\prec (N^*, \delta)_{\delta\in
 N\cap \nacc E}\}$. (only nacc replaces $S^{\lambda^+}_\lambda$).
\lz
\bigbreak
\centerline{\bf References}
\medskip
\item{[CK]} C.~C.~Chang, J. H. Keisler: Model Theory,
North-Holland, 1973.

\item{[M]} A.~Mekler: Universal structures in power $\aleph_1$,
{\sl Journal of Symbolic Logic \bf 55} (1990), 466--477.

\item{[Sh-c]}  S. Shelah :{\it Classification Theory and the number of 
non isomorphic models}, North Holland Publ.\ Co.,Studies in Logic and 
the 
Foundation of Math vol. 92, 1990, 705 + xxxiv.

\item{[Sh-g]} {\it Cardinal Arithmetic, }, Oxford University Press , 
to appear. 


\item{[Sh 80]}  A weak generalization of MA to higher
cardinals, {\it Israel J.of Math.}, {\bf 30} (1978), 297--306.

\item{[Sh 92]}  Remarks on Boolean algebras ,
{\it Algebra Universalis } 
{\bf 11} (1980), 77-89.

\item{[Sh 93]} Simple Unstable Theories, {\it Annals of Math.
Logic.},
 {\bf 19} (1980) 177--204.

\item{[Sh 100]} Independence results, {\it J. of Symb. Logic},
{\bf  45} (1980) 563--573.

\item{[Sh 107]} Models with second order properties IV, 
A general method and eliminating 
diamonds, {\it Annals. of Math. Logic}, {\bf 25} (1983) 183--212.

\item{[KfSh 150]}
M. Kaufman and S. Shelah: The Hanf number of stationary logic,
 {\it Notre Dame J. of Formal Logic}, {\bf 27} (1986) 111--
123.

\item{[ShLH 162]} S. Shelah,  C. Laflamne and B. Hart:
Models with second order properties V: A General principle,
{\it Annals of Pure and Applied Logic}, to appear.

\item{[GrSh 174]}
R. Grossberg and S. Shelah: On universal locally finite groups,
{\it Israel J. of Math.}, {\bf 44} (1983), 289--302.

\item{[Sh 175]}
On Universal graphs without instances of CH.
{\it Annals of Pure and Applied Logic}
,{\bf 26} (1984) 75--87.

\item{[Sh 175a]}
Universal graphs without instances of CH,
{\it Israel J. of Math.}, {\bf 70} (1990), 69--81.

%
%
\item{[Sh 288]}
Was Sierpenski right II, Proceedings of the 1/91 conference in 
Budapest in honour of A. Hajnal and V. T. S\'os,  Colloquia Mathemtica 
Societatis Janos 
Bolyai 60. Sets Graphs and Numbers , Budapest (Hungary )  ,1991,
637-668.

\item{[Sh 326]} Vive la Difference I,
{\it Proceedings of the Conference in Set Theory}, eds. H. Judah, W.
Just, and H. Woodin,   Springer Verlag
MSRI publications {\bf 26} (1992) , 357-405

\item{[Sh 365]} There are Jonsson algebras in many
inaccessible cardinals, Ch. IV of [Sh g].

\item{[Sh 405]}  Vive la difference II - refuting Kim's  
conjecture , {\it Israel J. of Mathematics } , accepted

\item{[KjSh 409]} M. Kojman and S. Shelah: 
Nonexistence of universal orders in many cardinals, {\it
Journal of Symbolic Logic \bf 57} (1992), 875--891.

\item{[Sh 413]} More Jonsson algebras and colorings, 
{\it Archive f.~Math.~Logic. } , accepted

\item {[KjSh 447]}  M. Kojman and S. Shelah: The universality
spectrum of stable unsuperstable theories, accepted to {\it Annals of
Pure and Applied Logic}. 
{\bf 58} (1992) , 57-92

\item{[KjSh 455]} M.Kojman and S. Shelah: Universal Abelian Groups.
{\it  `Israel J. of Mathematics} , submitted

\item{[KjSh 456]} in preparation

\item {[Sh 460]} The generalized continuum hypothesis revisited, 
preprint

\item{[Sh500]} Classifying Unstable Theories,  preprint.

\shlhetal
\bye